\documentclass[american]{amsart}
\usepackage[T1]{fontenc}
\usepackage[latin1]{inputenc}
\usepackage{hyperref}
\makeatletter

\theoremstyle{plain}    
\newtheorem{thm}{Theorem}[section]
\numberwithin{equation}{section} 
\numberwithin{figure}{section} 
\theoremstyle{plain}    
\newtheorem{cor}[thm]{Corollary} 
\newtheorem{lem}[thm]{Lemma} 
\theoremstyle{plain} 
\newtheorem{LKlem}[thm]{Local Key Lemma} 
\theoremstyle{plain}    
\newtheorem{prop}[thm]{Proposition} 
\theoremstyle{remark}

\theoremstyle{remark}    
 
\theoremstyle{plain}    
\newtheorem{Def}[thm]{Definition} 

\usepackage[all]{xy}
\xyoption{dvips}
\CompileMatrices

\makeatother

\begin{document}

\title{Numerically trivial foliations}


\author{Thomas Eckl}




\address{Thomas Eckl, Institut für Mathematik, Universität
  Bayreuth, 95440 Bayreuth, Germany}

\email{thomas.eckl@uni-bayreuth.de}

\urladdr{http://btm8x5.mat.uni-bayreuth.de/\~{}eckl}

\begin{abstract}
Given a positive singular hermitian metric of a pseudoeffective line bundle on
a complex K\"ahler manifold, a singular foliation is
constructed  satisfying certain analytic analogues of numerical conditions. 
This foliation refines Tsuji's numerically trivial fibration and the 
Iitaka fibration. Using almost positive singular hermitian metrics
with analytic singularities on a pseudo-effective line bundle , a foliation is
constructed refining the nef fibration. If the singularities of the foliation 
are isolated points, the codimension of the leaves is an upper bound to the
numerical 
dimension of the line bundle, and the foliation can be interpreted as a 
geometric reason for the deviation of nef and Kodaira-Iitaka dimension. 
Several surface examples are studied in more details, 
$\mathbb{P}^2$
blown up in 9 points giving a 
counter example to equality of numerical dimension and codimension of the 
leaves.
\end{abstract}

\maketitle
\tableofcontents

\bibliographystyle{alpha}


\section{Introduction}

\noindent
In the last few years several fibrations related to a nef or even 
pseudoeffective line bundle 
$L$
on a projective complex manifold were constructed 
whose fibers satisfy certain numerical properties with respect to a sometimes 
modified intersection theory: 

\noindent
For a nef line bundle
$L$,
the usual intersection theory is taken by \cite{BCEKPRSW00} to
define (and construct) the so called \textbf{nef fibration} whose fibers 
contain only curves 
$C$
with
$L.C = 0$.
The base dimension of this fibration is called the nef dimension of
$L$ , 
and it can be proven that it is never smaller than the numerical dimension
$\nu(L)$
of
$L$.
Note however that already for surfaces there are explicit counter examples to 
equality, cf.~section~\ref{MumEx-ssec}. 

\noindent
Even earlier, Tsuji~\cite{Tsu00} associated an intersection theory to positive
singular hermitian metrics 
$h$
on pseudoeffective line bundles 
$L$
by defining
\[ (L,h).C := \limsup_{m \rightarrow \infty} \frac{1}{m} 
              h^0(\tilde{C}, \mathcal{O}_{\tilde{C}}(m\pi^\ast L) \otimes 
                             \mathcal{I}((\pi^\ast h)^m)). \]
Here, 
$\pi: \tilde{C} \rightarrow C$
is the normalization of an irreducible curve
$C$
not contained in the singular locus of
$h$,
and
$\mathcal{I}((\pi^\ast h)^m))$
denotes the multiplier ideal sheaf of the pulled back metric
$(\pi^\ast h)^m$
on
$\tilde{C}$.
A projective complex manifold is called numerically trivial in Tsuji's sense
iff 
$(L,h).C = 0$
for all such curves 
$C$.
In~\cite{E02} other possible definitions of these intersection numbers are 
discussed, their relations are studied, and the fibration map with numerically 
trivial fibers is constructed, according to the suggestions of Tsuji. 
  
\noindent
Finally, Takayama~\cite{Tak02} defined intersection numbers reflecting 
properties of the linear sytems
$|mL|$
by using the asymptotic multiplier ideal sheaf
$\mathcal{J}(||L||)$:
\[ ||L,C|| := \lim_{m \rightarrow \infty} m^{-1} \deg_C mL \otimes 
                                                 \mathcal{J}(||mL||), \] 
where
$C$
is an irreducible curve not contained in the stable base locus 
$\bigcap_{m \in \mathbb{N}} \mathrm{Bs}|mL|$
of
$L$.
The resulting fibration turns out to be the well known Iitaka fibration.

\noindent
The motivation for this work is to give a more unified treatment of all these 
fibrations and to give geometric reasons for the deviation of nef, numerical, 
and Kodaira-Iitaka dimension of a nef line bundle on a projective manifold. 
Three surface examples will illustrate the ideas developed to this purpose.

\noindent
The first example is due to Mumford and has the property that the nef 
dimension is bigger than the numerical dimension: Start with a 
smooth projective curve
$C$
of genus
$\geq 2$
with the unit circle 
$\Delta$
as universal covering and an irreducible unitary representation
$\rho: \pi_1(C) \rightarrow GL(2,\mathbb{C})$
of the fundamental group of
$C$.
This defines a rank 2 vector bundle 
$E = (\Delta \times \mathbb{C}^2)/\pi_1(C)$
on
$C$
of degree 0 where the action of
$\pi_1(C)$
is given by covering transformations on 
$\Delta$
and the representation
$\rho$
on
$\mathbb{C}^2$. 

\noindent
Mumford proved that the nef line bundle 
$L = \mathcal{O}_{\mathbb{P}(E)}(1)$ 
on the projectivized bundle 
$\mathbb{P}(E)$
is stable hence the restriction of 
$L$ 
to all curves 
$D \subset \mathbb{P}(E)$
is positive. On the other hand 
$\deg E = 0$
hence 
$L.L = 0$. 
Hence the numerical dimension
$\nu(L)$
is
$1$,
while the nef reduction map is the identity, and the nef dimension is
$2$.

\noindent
It seems quite obvious how to explain this deviation: the ruled surface 
$\mathbb{P}(E)$
carries a foliation induced by the images of the
$\Delta \times l$
in
$\mathbb{P}(E)$
(where
$l$
is a line through the origin in
$\mathbb{C}^2$).
Furthermore, locally the leaves of this foliation are mapped to points by the 
morphism induced by
$|L|$,
which is a kind of numerical triviality. 

\noindent
This motivates the construction of numerically trivial foliations w.r.t. some
positive closed current (which may be the curvature current of some hemitian 
metric on a nef line bundle) on a complex manifold
$X$.
The starting point is an interesting criterion for numerical triviality in
Tsuji's sense:
\begin{thm} \label{NumTriv-thm}
Let
$X$
be a smooth projective complex manifold, let
$L$
be a pseudo-effective line bundle on
$X$
with positive singular hermitian metric 
$h$
such that
$X$
is
$(L,h)-\!$
numerically trivial. Then the curvature current
$\Theta_h$
may be decomposed as
\[ \Theta_h = \sum_i a_i [D_i] \]
where the
$D_i$
form a countable set of prime divisors on
$X$ 
and the 
$a_i$
are
$> 0$. 
\end{thm}

\noindent
This is proven in~\cite{E02}, and by trivial arguments the converse of this 
theorem is also true. It shows that 
numerical triviality is a local property of currents and does not depend on 
projectivity. Hence it is possible to localize the notion of 
numerically trivial fibrations to the notion of a foliation
with numerically trivial leaves (details in 
sections~\ref{NumTriv-ssec},\ref{MaxNTFol-ssec}). 
The main result is the following:
\begin{thm} \label{NTFol-thm}
On a (not necessarily compact) complex manifold 
$X$
with a positive closed
$(1,1)-\!$
current
$T$
there exists a 
\textbf{maximal} foliation with numerically trivial leaves w.r.t.
$T$, 
that is
the leaves of every foliation with numerically trivial leaves are contained in 
leaves of this foliation. 
\end{thm}

\noindent
It is called the numerically trivial foliation 
w.r.t.
$T$.
The construction rests essentially on the Local Key 
Lemma which allows to unite different foliations with numerically trivial 
leaves, and the proof of this lemma is an easy consequence of another 
interpretation of numerically trivial fibrations 
$f: X \rightarrow Y$
w.r.t. to some closed positive
$(1,1)-\!$
current
$T$:
The residue current 
$R$
of the Siu decomposition
$T = \sum a_i [D_i] + R$
must be the pull back of a (positive) current on
$Y$ 
(details in~\ref{MaxNTFol-ssec}).

\noindent
If
$X$
is projective and
$T$
the curvature current of a positive singular hermitian metric on a line bundle,
Tsuji's numerically trivial fibration will be the fibration maximal among 
those whose fibers are contained in the leaves of the numerically trivial 
foliation w.r.t. 
$T$;
details in section~\ref{NTFib-ssec}.
The same construction gives the Iitaka fibration of a line 
bundle 
$L$
with Kodaira-Iitaka dimension
$\kappa(L) \geq 0$
provided one uses the positive singular hermitian metric
$h$
on
$L$
defined as
\[ h = \limsup_{m \rightarrow \infty} (h_{|mL|})^{\frac{1}{m}} \]
where 
$h_{|mL|}$
is the (singular) hermitian metric on 
$|mL|$
defined by the global sections of
$mL$ 
(see~\cite{Tsu99}).
In this case even more is true: The numerical trivial foliation w.r.t. 
$h$
is already the Iitaka fibration (section~\ref{ITFib-ssec}). 

\noindent
It is not possible to find a positive singular hermitian metric which defines
the nef fibration in a similar way, as shown by an example constructed 
in~\cite{DPS94} which is quite similar to Mumford's example 
(section~\ref{firstex-ssec}): Start with an elliptic curve 
$C$
and take as the rank 2 vector bundle 
$E$
the unique nontrivial extension of the structure sheaf
$\mathcal{O}_C$.
As in Mumford's example the numerical dimension of the nef line bundle
$L = \mathcal{O}_{\mathbb{P}(E)}(1)$
is 
$1$,
while the nef dimension is
$2$.

\noindent
The remarkable feature of this example is the fact that the only positive 
singular hermitian metric on
$L$
is given by the unique section of
$L$, 
the ``section at infinity'' of
$\mathbb{P}(E)$
(proof in~\cite{DPS94}), and
$\mathbb{P}(E)$
is numerically trivial w.r.t. this metric. Hence the obvious foliation on
$\mathbb{P}(E)$
induced by the universal cover
$\mathbb{C}$
of
$C$
cannot be interpreted as the numerically trivial foliation w.r.t. some positive
metric on
$L = \mathcal{O}_{\mathbb{P}(E)}(1)$.

\noindent
Ideas how to deal with this situation 
may be found in Boucksom's construction of a divisorial 
Zariski decomposition and his definition of ``moving'' intersection numbers 
on pseudoeffective line bundles~\cite{Bou02} on compact K\"ahler manifolds. 
Both notions show that it is 
extremely useful to loosen the restriction on positivity and to 
consider sequences of almost positive 
$(1,1)-\!$
currents in a fixed cohomology class
$\alpha$
whose negative parts tend to 0. 

\noindent
For nef line bundles the moving intersection numbers coincide
with the usual ones. In particular, if
$C \subset X$
is a smooth compact curve on a compact K\"ahler manifold
$X$
with K\"ahler form
$\omega$,
and 
$L$
is a nef line bundle with first Chern class
$\alpha := c_1(L) \in H^{1,1}(X,\mathbb{R})$,
\[ L.C = \lim_{\epsilon \downarrow 0} \sup_T \int_{X - \mathrm{Sing} T}
         (T + \epsilon \omega) \wedge [C]  
       = \lim_{\epsilon \downarrow 0} \sup_T \int_{C - \mathrm{Sing} T}
         (T + \epsilon \omega), \]
where the 
$T$'s
run through all closed currents representing
$\alpha$
such that
$T \geq - \epsilon \omega$,
and
$[C]$
is the integration current of the submanifold
$C$
of bidegree
$(n-1,n-1)$.
(For further details see section~\ref{MovIntNum-ssec}.) 

\noindent
It is obvious that in this case
$L.C = 0$
iff
$\lim_{\epsilon \downarrow 0} \sup_T \int_{\Delta - \mathrm{Sing} T}
         (T + \epsilon \omega) = 0$
for all disks
$\Delta \subset C$.
Thus it is justified to interpret numerical triviality w.r.t. a 
pseudo-effective class as a local property: An immersed submanifold
$Y \subset X$
(closed or not) is called numerically trivial w.r.t. a pseudo-effective class
$\alpha$
iff
\[ \lim_{\epsilon \downarrow 0} \sup_T \int_{\Delta - \mathrm{Sing} T}
         (T + \epsilon \omega) = 0 \]
(where the 
$T$'s
run through all closed currents on 
$X$
contained in
$\alpha[-\epsilon\omega]$)
for all holomorphically immersed disks
$\Delta \subset Y$.
And a foliation will be called numerically trivial w.r.t. 
$\alpha$
iff (locally) almost every  leaf is numerically trivial w.r.t. 
$\alpha$.

\noindent
It is possible to prove an analog to the Local Key Lemma, hence there is a
maximal numerically trivial foliation w.r.t.
$\alpha$.
It is contained in every numerically trivial foliation w.r.t. a positive 
current representing
$\alpha$. 
If
$\alpha$
is the first Chern class of some nef line bundle
$L$
on a projective manifold
$X$,
the nef fibration of 
$L$
is the maximal fibration contained in the foliation (which will be called the 
\textbf{nef foliation} in that case). Furthermore, the Kodaira-Iitaka 
fibration contains the nef foliation, and one 
gets a nice geometric reason for deviations of the Kodaira-Iitaka 
and the nef dimension of nef line bundles on projective manifolds:
$\kappa(L) < n(L)$
if the nef foliation is not a fibration. It is a very interesting open 
question whether the converse of this statement 
is also true. More generally: Is the fibration with the smallest fiber 
dimension which contains the nef foliation the Kodaira-Iitaka fibration ?

\noindent
Finally, it is shown that the codimension of the leafs is an upper bound for
the numerical dimension of
$\alpha$, 
if the singularities of the foliation are isolated points. It is not clear to 
the author how to weaken this assumption or if there are counter examples. To
get better answers it seems necessary to have a closer look at the structure
of numerical trivial foliations around the singularities.

\noindent
The last section of the paper constructs nef foliations
of nef line bundles on surfaces. The first two
examples are those due to Mumford and Demailly-Peternell-Schneider. 
In Mumford's example it is easy to construct a smooth 
closed positive
$(1,1)-\!$
current on
$L = \mathcal{O}_{\mathbb{P}^1}(1)$
such that the associated nef foliation is the obvious one: 
Take a measure
$\omega$
invariant w.r.t. the representation of
$\pi(C)$
in
$\mathrm{PGL}(2)$.
This gives a measure on
$(\Delta \times \mathbb{P}^1)/\pi(C)$
transversal to the foliation induced by the images of
$\Delta \times \left\{ p \right\}$. 
Averaging out the integration currents of the leaves with this transverse 
measure gives an (even smooth) closed positive
$(1,1)-\!$
current in the first Chern class of
$L = \mathcal{O}_{\mathbb{P}(E)}(1)$
which vanishes on the leaves but not in any transverse direction.

\noindent
The Demailly-Peternel-Schneider example is more difficult: A complicated 
glueing argument leads to almost positive currents which determine the obvious
foliation. 

\noindent
The last example deals with 
$\mathbb{P}^2$
blown up in
$9$
points and is interesting in many ways. In particular, if one fixes 8
points in sufficiently general position, varying the last point will give a 
nef fibration in the torsion points, but there is no nef foliation on the 
whole family with 
$1$-
dimensional leaves. Hence, the nef fibrations in varieties over torsion points
do not converge against a foliation in varieties over (general) non-torsion 
points. This somehow answers a question asked in \cite{DPS96}.   

\subsubsection*{Acknowledgement }

This article was mainly written during two stays at the Institut Fourier in 
Grenoble. The first was paid by the DFG-Schwerpunkt "Global Methods in Complex
Geometry", the second by the Institut Universitaire de France. The author 
enormously benefitted from many discussions with J.-P. 
Demailly and S. Boucksom who generously shared their new insights on moving 
intersection numbers with him and gave the author a lot of encouragement.
Furthermore, the author enjoyed a lot the warm hospitality created by the whole
institute.

\section{Numerically trivial foliations} \label{NTFol-sec}

\subsection{Numerical triviality} \label{NumTriv-ssec}

\noindent
As proposed in the introduction numerical triviality of (not necessarily 
compact) complex manifolds is defined via the criterion of 
theorem~\ref{NumTriv-thm}:
\begin{Def}
Let
$X$
be a complex manifold and 
$T$
a positive closed
$(1,1)-\!$
current on
$X$.
Then
$X$
is called numerically trivial w.r.t. 
$T$
iff 
\[ T = \sum a_i[D_i], \]
for countably many prime divisors 
$D_i$
in
$X$
and real numbers
$a_i \geq 0$.
\end{Def}

\noindent
To compare later on Tsuji's numerically trivial fibration with the numerically
trivial foliation it is useful to define numerical triviality for any 
irreducible analytic subsets (not only for submanifolds):
\begin{Def} \label{NumTrivSubVar-def}
Let
$X$
be compact complex manifold and
$\Theta$
a positive closed
$(1,1)-\!$
current on
$X$.
Let
$Y \subset X$
be a positive dimensional analytic subset of
$X$ 
such that
$\Theta$
may be restricted to
$Y_{reg}$,
the smooth part of
$Y$.
Then
$Y$
is called numerically trivial with respect to
$\Theta$
iff for all holomorphic maps
$f: \Delta^k \rightarrow Y$
such that
$f^\ast \Theta$
exists the complex manifold
$\Delta^k$
is numerically trivial with respect to
$f^\ast \Theta$.
\end{Def}

\noindent
This definition is consistent with the definition of numerically 
trivial manifolds:
\begin{prop}
Let
$X$
be a complex manifold and 
$\Theta$
a positive closed
$(1,1)-\!$
current.
$X$
is numerically trivial w.r.t.
$\Theta$
iff for all holomorphic maps
$f: \Delta^k \rightarrow X$
such that
$f^\ast \Theta$
exists 
$\Delta^k$
is numerically trivial w.r.t. 
$f^\ast \Theta$. 
\end{prop}
\begin{proof}
The ``only if'' part is a trivial consequence of the equality
$f^\ast([D_i]) = [f^\ast(D_i)]$
for (integration currents of) divisors. The other direction follows from the 
Siu-decomposition \cite[(2.18)]{Dem00}
\[ \Theta = \sum_i \nu_i[D_i] + R, \]
where
$R$
is a positive closed
$(1,1)-\!$
current such that the Lelong number level sets
$E_c(R)$
have no codim 1 components. If
$R \not\equiv 0$
there will exist an open set
$U \cong \Delta^n$
in
$X$
such that
$R_{|U} \not\equiv 0$,
hence
$U$
is not numerically trivial w.r.t.
$\Theta_{|U}$.
\end{proof} 

\noindent
The definition of numerical triviality can be further simplified by means of 
the following proposition:
\begin{prop} \label{NumTrivCrit-prop}
Let
$X$
be a complex manifold and 
$\Theta$
a positive closed
$(1,1)-\!$
current.
$X$
is numerically trivial w.r.t.
$\Theta$
iff for all holomorphic maps
$f: \Delta \rightarrow X$
such that
$f^\ast \Theta$
exists 
$\Delta$
is numerically trivial w.r.t. 
$f^\ast \Theta$. 
\end{prop}
\begin{proof}
The ``only if'' part follows by definition. For the other direction start again
with the Siu decomposition
$\Theta = \sum_i \nu_i[D_i] + R$.
Let
$\Delta^n \cong U \subset X$
be an open subset and let
$q: \Delta^n \rightarrow \Delta^{n-1}$
be the projection onto the first
$n-1$
factors. Since the Lelong number level sets
$E_c(R)$
contain no codim 1 component, very general fibers 
$F$
of 
$q$
do not intersect any of the 
$E_c(R)$.
By the results of \cite{BMElM00} there is a pluripolar set
$N \subset \Delta^{n-1}$
such that the level sets
$E_c(R_{|F}) = \emptyset$
for the restriction of 
$R$
to all fibers 
$F$
over points outside of
$N$.
By assumption 
$R_{|F} \cong 0$.

\noindent
By the following lemma there exists a positive closed
$(1,1)-\!$
current
$S$
on
$\Delta^{n-1}$
such that
$R = q^\ast S$.
Let
$D = \Delta^{n-1} \times \left\{ p \right\}$
be a section of
$q$
such that
$R_{|D}$
is well defined. By induction
$R_{|D} \equiv 0$.
Since the projection
$q: D \rightarrow \Delta^{n-1}$
is an isomorphism 
$S \equiv 0$
hence
$R \equiv 0$.
\end{proof}

\begin{lem} \label{pullback-lem}
Let
$T$
be a positive closed 
$(1,1)-\!$
current on
$\Delta^n$ 
and let
$q: \Delta^n \rightarrow \Delta^{n-1}$
be the projection onto all factors but the last one. If
$T_{|q^{-1}(x)} \equiv 0$
for all
$x$
outside a pluripolar set
$N \subset \Delta^{n-1}$
then there will be a positive closed 
$(1,1)-\!$
current 
$S$
on
$\Delta^{n-1}$
such that
$T = q^\ast S$.
\end{lem}
\begin{proof}
The positive current
$T$
may be written as
\[ T = i \sum_{i,j} \Theta_{ij} dz_i \wedge d\overline{z}_j \]
where the 
$\Theta_{ij}$
are complex measures on
$\Delta^n$
(\cite[(1.15)]{Dem00}). That
$T$
is a real current implies 
$\Theta_{ij} = \overline{\Theta}_{ji}$. 
Since
$T$
is positive,
$\sum \lambda_i \overline{\lambda_j} \Theta_{ij}$
is a positive measure for all vectors
$(\lambda_1, \ldots, \lambda_n) \in \mathbb{C}^n$.
Hence
\[ \lambda_{i}\overline{\lambda_i} \Theta_{ii} +
   \lambda_{i}\overline{\lambda_n} \Theta_{in} +
   \lambda_{n}\overline{\lambda_i} \Theta_{ni} +
   \lambda_{n}\overline{\lambda_n} \Theta_{nn} \geq 0\ \ 
   \forall (\lambda_i,\lambda_n) \in \mathbb{C}^2. \]

\noindent
\textit{Claim.} As a
$(1,1)-\!$
current
$i \Theta_{nn}dz_n \wedge d\overline{z}_n = 0$.
\begin{proof}
By definition one has to show that
\[ \int_{\Delta^n} i \Theta_{nn}dz_n \wedge d\overline{z}_n \wedge \alpha
                               idz_1 \wedge d\overline{z}_1 \wedge \ldots
                        \wedge idz_{n-1} \wedge d\overline{z}_{n-1} = 0 \]
for all complex valued functions
$\alpha \in \mathcal{C}^\infty_c(\Delta^n)$. 
Since
$T_{|q^{-1}(x)} = i \Theta_{nn}dz_n \wedge d\overline{z}_n$
\begin{eqnarray*} 
  \lefteqn{\int_{\Delta^n} i \Theta_{nn}dz_n \wedge d\overline{z}_n \wedge 
                           \alpha
                               idz_1 \wedge d\overline{z}_1 \wedge \ldots
                        \wedge idz_{n-1} \wedge d\overline{z}_{n-1} =}\\
   & & \ \ \ \ \int_{\Delta^n} T \wedge \alpha
                               idz_1 \wedge d\overline{z}_1 \wedge \ldots
                        \wedge idz_{n-1} \wedge d\overline{z}_{n-1}, 
\end{eqnarray*}
and the slicing formula \cite[(1.22)]{Dem00} implies that this is equal to
\[ \int_{\Delta^{n-1}} \left( \int_{q^{-1}(x^\prime)} T_{|q^{-1}(x^\prime)} 
                       \wedge \alpha_{|q^{-1}(x^\prime)} \right)
                       idz_1 \wedge d\overline{z}_1 \wedge \ldots
                        \wedge idz_{n-1} \wedge d\overline{z}_{n-1}. \]
This is
$0$
because
$T_{|q^{-1}(x)} \equiv 0$
for all
$x$
outside a pluripolar set
$N \subset \Delta^{n-1}$.
\end{proof}

\noindent
Consequently,
\[ \Theta_{ii} + \overline{\lambda}_n \Theta_{in} + \lambda_n \Theta_{ni} =
   \Theta_{ii} + \overline{\lambda}_n \overline{\Theta}_{ni} + 
                 \lambda_n \Theta_{ni}  
   \geq 0 \]
for all
$\lambda_n \in \mathbb{C}$.
Now suppose that
$\Theta_{ni} \neq 0$, 
i.e. there is a smooth real valued function
$\alpha \geq 0$
with compact support such that
$\Theta_{ni}(\alpha) \neq 0$.
Then there is a 
$\lambda_n \in \mathbb{C}$
such that
\[ \Theta_{ii}(\alpha) + 
   \overline{\lambda}_n \overline{\Theta_{ni}(\alpha)} + 
   \lambda_n \Theta_{ni}(\alpha)  < 0. \]

\noindent
This is a contradiction. Hence
$\Theta_{in} = \Theta_{ni} = 0$
for all
$i \leq n-1$.

\noindent
Next, the closedness of 
$T$
implies 
\[ \frac{\partial}{\partial z_n} \Theta_{ij} = 
   \frac{\partial}{\partial \overline{z_n}} \Theta_{ij} = 0 \ \ \ 
   \forall i,j \leq n-1. \]
Hence the
$\Theta_{ij}$
only depend on
$z_1, \ldots, z_{n-1}$.
One finally gets
\[ T = q^\ast S = i \sum_{i,j \leq n-1} \Theta_{ij} dz_i \wedge 
                                          d\overline{z_j} \]
and
$S$
is a closed positive
$(1,1)-\!$
current on
$\Delta^{n-1}$. 
\end{proof}

\noindent
Proposition~\ref{NumTrivCrit-prop} has an easy
\begin{cor} \label{numtriv-cor}
Let
$X$
be a complex manifold and 
$\Theta$
a positive closed 
$(1,1)-\!$
current.
Let
$Y \subset X$
be an irreducible analytic subset such that
$\Theta$
may be restricted to
$Y_{reg}$,
the smooth part of
$Y$.
Then
$Y$
is numerically trivial w.r.t.
$\Theta$
iff for all holomorphic maps
$f: \Delta \rightarrow Y \subset X$
such that
$f^\ast \Theta$
exists the complex manifold
$\Delta$
is numerically trivial w.r.t. 
$f^\ast \Theta$. \hfill $\Box$ 
\end{cor}

\noindent
As a consequence one can give an alternative definition of numerically trivial
irreducible analytic subsets using embedded resolutions of singularities by 
blowing up smooth centers. Such resolutions exist for arbitrary complex 
manifolds, at least on relatively compact open subsets 
(\cite{Hir64},\cite{BM97}).
\begin{prop}
Let
$X$
be a complex manifold and 
$\Theta$
a positive closed
$(1,1)-\!$
current. Let
$Y \subset X$
be an irreducible analytic subset such that
$\Theta$
may be restricted to
$Y_{reg}$,
the smooth part of
$Y$.
Then
$Y$
is numerically trivial w.r.t.
$\Theta$
iff for an embedded resolution
$f: \widetilde{Y} \rightarrow Y$
the complex manifold
$\widetilde{Y}$
is numerically trivial w.r.t. the pulled back current
$f^\ast \Theta$. 
\end{prop}
\begin{proof}
By the universal property of the blowup a map
$f: \Delta \rightarrow Y$
will factorize through the blow up 
$\pi: \widetilde{Y} \rightarrow Y$
of a smooth center if its image is not 
contained in the center. Furthermore the exceptional divisor is a 
projectivized bundle hence locally trivial. So at least locally  
there will be a map
$\widetilde{f}: \Delta \rightarrow \widetilde{Y}$
such that
$\pi \circ \widetilde{f} = f$
if the image is contained in the center.   
\end{proof}

\noindent
Finally there is a useful criterion for numerical triviality:
\begin{prop} \label{gennumtriv-prop}
Let
$X$
be a complex manifold and let
$\Theta$
be an almost positive
$(1,1)-\!$
current. Then
$X$
is numerically trivial w.r.t.
$\Theta$
iff there is an analytic subset
$A \subset X$
such that
$X-A$
is numerically trivial w.r.t.
$\Theta$.
\end{prop}
\begin{proof}
This is a direct consequence of the following standard arguments. First a 
closed 
$(1,1)-\!$
current is 0 iff it is already 0 outside a set of real codimension 4. 
\cite[(1.21)]{Dem00}. Second, for complete pluripolar sets 
$E$
(as are analytic 
subsets) 
$\Theta = \mathbf{1}_{X-E} \Theta + \mathbf{1}_E \Theta$:
This is true for the closed positive current 
$\Theta + C\omega$
by \cite[(1.19)]{Dem00},
hence also for
$\Theta$. 
But for
$E$
a codimension 1 analytic subset
$\mathbf{1}_E \Theta = m_E [E]$
where
$m_E$
is the generic Lelong number on
$E$
\cite[(2.17)]{Dem00}.   
\end{proof}

\subsection{Existence of maximal numerically trivial foliations}
\label{MaxNTFol-ssec}

\noindent
Consider singular foliations as described in the Appendix:
\begin{Def}\label{ntfol-def}
Let
$X$
be an
$n-\!$
dimensional compact complex manifold with hermitian metric  
$\omega$
and
$\Theta \geq 0$
a positive  closed 
$(1,1)-\!$
current on
$X$.
A singular foliation
$\{\mathcal{F},(U_i,p_i)\}$
is said to induce a (singular) numerically 
trivial foliation w.r.t. 
$\Theta$
iff almost every fiber of  
$p_i$
is numerically trivial w.r.t.
$\Theta$.
\end{Def}

\noindent
Note that the condition about the fibers of the
$p_i$
is much stronger than in the case of numerical trivial 
fibrations: Here it was only necessary to assume that the union of numerically
trivial fibers is not a Lebesgue zero set, and the numerical triviality of all
fibers over points lying in the complement of a pluripolar set followed. In 
the foliation case there exist counter examples to this conclusion: On 
$\Delta^2$
consider the plurisubharmonic function
\[ \phi(z_1,z_2) = \max(\log(1+|z_1z_2|^2)-\log \frac{5}{4}, 0). \]
For every 
$z_2-\!$
fiber 
$F$
with
$|z_2| < \frac{1}{2}$
the restriction
$\phi_{|F}$
is
$\equiv 0$.
But for 
$|z_1z_2| > \frac{1}{2}$
one sees that
$\phi \equiv \log(1+|z_1z_2|^2)-\log \frac{5}{4}$.

\noindent
Theorem~\ref{NTFol-thm} states the existence of a maximal 
numerically trivial foliation
with respect to the inclusion relation
``$\sqsubset$'' of singular foliations, see the Appendix.
The strategy to prove the existence of this maximal foliation is 
essentially the 
same as for the existence proof of numerical trivial fibrations:
one proves that the common refinement 
$\left\{ \mathcal{H}, (W_k, r_k: W_k \rightarrow \Delta^{n-m}) \right\}$
of two numerically trivial foliations
$\left\{ \mathcal{F}, (U_i, p_i: U_i \rightarrow \Delta^{n-k}) \right\}$,
$\left\{ \mathcal{G}, (V_j, q_j: V_j \rightarrow \Delta^{n-l}) \right\}$,
(see the Appendix) is again a numerically trivial foliation.

\noindent
The main step is to establish a local analog to the Key 
Lemma  in~\cite{E02}. It is stated for the following configuration: Let
$W \cong \Delta^n$
be a complex manifold with two projections
$p_1: W \rightarrow \Delta^{n-k}$,
$p_2: W \rightarrow \Delta^{n-l}$
such that a smallest projection
$p: W \rightarrow \Delta^{n-m}$
as constructed in the Appendix exists.
\begin{LKlem} \label{LKlem}
If the fibrations induced by
$p_1$
and
$p_2$
are numerically trivial w.r.t. a positive closed
$(1,1)-\!$
current 
$\Theta$
on
$W$,
then the foliation induced by
$p$
will also be numerically trivial w.r.t. 
$\Theta$.
\end{LKlem}

\noindent
Since any two points on
$\Delta^{n-k}$
may be connected by a sequence of images of
$p_2-\!$
fibers this is a consequence of
\begin{lem}
Let
$\Theta$
be a positive closed
$(1,1)-\!$
current on
$\Delta^n$, 
let
$q: \Delta^n \rightarrow \Delta^k$
be the projection onto the last 
$k$
factors and let
\[ V = \{x \in \Delta^n|x_{1} = \ldots = x_l = 0\} \subset \Delta^n \]
be an analytic subset mapping surjectively on
$\Delta^k$
(that is, 
$l < k$). 
If almost every 
$q-\!$
fiber and 
$V$
are numerically trivial w.r.t. 
$\Theta$
then
$\Delta^n$
will be numerically trivial w.r.t. 
$\Theta$.
\end{lem} 
\begin{proof}
Let
$\Theta = \sum_i a_i[D_i] + R$
be the Siu decomposition. Lemma~\ref{pullback-lem} shows that
$R = q^\ast S$
for some closed positive 
$(1,1)-\!$
current 
$S$
on
$\Delta^k$.
Now, at least locally each map
$f: \Delta \rightarrow \Delta^k$
is liftable to a map 
$\widetilde{f}: \Delta \rightarrow V$,
that is
$f = q \circ \widetilde{f}$
(this is obvious for projections). Hence
$S$
is numerically trivial by the criterion in proposition~\ref{gennumtriv-prop}.
But divisorial components in the Siu-decomposition of
$S$
would give divisorial components of the Lelong number level sets of
$R = q^\ast S$. 
Therefore
$S \equiv 0$
hence
$R \equiv 0$.
\end{proof}

\noindent
Note again that the Local Key Lemma needs stronger assumptions on the
fibers than the Key Lemma. This is shown by the
same counter example as above: The horizontal sections 
$\{ z_2 = a \}$
are also numerically trivial as long as~$|a| < \frac{1}{2}$. 

\noindent
Now it is an easy consequence of the Local Key Lemma to show that common 
refinement 
$\left\{ \mathcal{H}, (W_k, r_k: W_k \rightarrow \Delta^{n-m}) \right\}$
of two numerically trivial foliations
$\left\{ \mathcal{F}, (U_i, p_i: U_i \rightarrow \Delta^{n-k}) \right\}$,
$\left\{ \mathcal{G}, (V_j, q_j: V_j \rightarrow \Delta^{n-l}) \right\}$
is again numerically trivial.

\noindent
This ends the proof of Theorem~\ref{NTFol-thm}.

\subsection{Tsuji's numerically trivial fibrations} \label{NTFib-ssec}

Now let
$X$
be a smooth projective complex manifold and
$L$
a pseudoeffective holomorphic line bundle on
$X$
with positive singular hermitian metric 
$h$. 
As already mentioned in the beginning, the notion of numerical triviality used
to construct Tsuji's
numerically trivial fibrations is derived from an intersection number
$(L,h).C$
of an irreducible curve
$C \subset X$
not contained in the singular locus of
$h$
with the pair
$(L,h)$.
A subvariety
$Y \subset X$
is numerically trivial iff 
$(L,h).C = 0$
for all irreducible curves
$C \subset X$.
The analysis of these intersection numbers in \cite{E02} shows that
\[ (L,h).C = (\pi^\ast L, \pi^\ast h).\overline{C} = 
   \pi^\ast L.\overline{C} - \sum_{x \in \overline{C}} \nu(\pi^\ast h,x), \]
where
$\pi: \overline{C} \rightarrow C$
is the normalization. In particular, if
$(L,h).C = 0$,
the curvature current of
$\pi^\ast h$
on
$\overline{C}$
may be written as
$\sum_{x \in \overline{C}} \nu(\pi^\ast h,x)[x]$.
Hence proposition~\ref{NumTrivCrit-prop} shows that numerically trivial 
subvarieties in the sense of definition~\ref{NumTrivSubVar-def}
are also numerically trivial in the sense just described.

\noindent
The converse is also true: By the birational invariance of numerical 
triviality \cite[2.6]{E02} the normalization and desingularization
$\overline{Y}$
of a numerically trivial subvariety
$Y$
is also numerically trivial (in Tsuji's sense). Hence the curvature current of
the pulled back metric is of the form 
$\sum \nu_i [D_i]$,
by theorem~\ref{NumTriv-thm}. But this implies certainly numerical triviality
of
$\overline{Y}$
in the sense of definition~\ref{NumTrivSubVar-def}, and since every holomorphic
map
$f: \Delta \rightarrow Y$
may be lifted to a holomorphic map
$f: \Delta \rightarrow \overline{Y}$,
the numerical triviality of
$Y$
follows.

\noindent
Now remind the construction of Tsuji's numerically trivial
fibration: It is the (up to birational equivalence unique) element with maximal
fiber dimension in the set of families
$\widetilde{f}: \mathfrak{X} \rightarrow \mathcal{N}$
with the following properties (\cite[3.3]{E02}):
\begin{itemize}
\item[(i)]
$\mathfrak{X} \subset X \times \mathcal{N}$,
$\mathfrak{X}, \mathcal{N}$
quasi-projective, irreducible, general fibres are subvarieties of
$X$;

\smallskip

\item[(ii)]
the projection
$p: \mathfrak{X} \rightarrow X$
is generically finite;

\smallskip

\item[(iii)]
$(L,h)$
is defined and not numerically trivial on sufficiently general fibres of 
$\widetilde{f}$,
i.e. on a set of fibres
$\mathcal{M} \subset \mathcal{N}$
which has not Lebesgue measure 
$0$;

\smallskip

\item[(iv)]
the fibres are generically unique, i.e. if
$U \subset \mathcal{N}$
is an open subset such that
$\widetilde{f}_{|U}$
is flat then the induced map 
$U \rightarrow \mathrm{Hilb}(X)$
will be generically bijective.
\end{itemize}

\noindent
It is shown that for the maximal element, the projection
$p: \mathfrak{X} \rightarrow X$
is really birational, and that all fibers where
$h_{|F} \not\equiv \infty$
are numerically trivial. But by the observation above, such a fibration can be
interpreted as a numerically trivial foliation
$\{\mathcal{F},(U_i,p_i)\}$:
Take
$\mathcal{F}$
as 
$p_\ast T_{\mathfrak{X}/\mathcal{N}}$, 
and let
$Z \subset X$
be an algebraic subset of points where
$p$
is an isomorphism and
$\widetilde{f}$
is smooth. Then
$X-Z$
may be covered by (analytically) open sets
$U_i$
such that there exist maps
$p_i: U_i \rightarrow \Delta^{n-k}$
with
$T_{U_i/\Delta^{n-k}} = \mathcal{F}_{|U_i}$.
This implies that the
$p_i-\!$
fibers are numerically trivial. Consequently, it is possible to characterize
Tsuji's numerically trivial fibration in the following way:
\begin{prop}
Let
$X$
be a smooth projective complex manifold and
$L$
a pseudoeffective holomorphic line bundle on
$X$
with positive singular hermitian metric 
$h$.
Then the birational fibration with maximal fiber dimension contained in the 
numerically trivial foliation w.r.t. the curvature current
$\Theta_h$
is Tsuji's numerically trivial fibration w.r.t. 
$(L,h)$.  \hfill $\Box$
\end{prop}

\subsection{The Iitaka fibration} \label{ITFib-ssec}

Let
$X$
be a projective complex manifold and
$L$
a line bundle with non-negative Kodaira-Iitaka dimension
$\kappa(X,L) \geq 0$.
Consider the set 
$N(L)$
of all
$m \in \mathbb{N}$
such that the linear systems
$|mL| \neq \emptyset$. 
Let
$m_0$
be the greatest common divisor of the numbers in
$N(L)$.
Then there is a positive integer
$m(L)$
such that
$|mm_0L| \neq \emptyset$
for all positive integers 
$m \geq m(L)$.
Choose generating sets 
$f_1, \ldots, f_{k_m}$
for the linear systems
$|mm_0L| \neq \emptyset$
and let
$h_m$
be the (possibly singular) hermitian metric on
$L$
with plurisubharmonic weight (on the base 
$\Omega \subset \mathbb{C}^n$
of a local trivialization 
$L \cong \Omega \times \mathbb{C}$)
\[ \phi_m = \frac{1}{2mm_0} \log(\sum_{i=1}^{k_m} |f_i|^2) \]
and curvature current
$\Theta_m = i\partial\overline{\partial}\phi_m$
(on
$\Omega$).
Let
$h_L$
be a smooth hermitian metric on
$L$
with weight 
$\phi_L$
on
$\Omega$
and smooth curvature form
$\Theta_L$.
Write
$\Theta_m = \Theta_L + i\partial\overline{\partial}\phi_m^\prime$
and normalize the 
$\phi_m^\prime$
by subtracting (if necessary) a positive constant 
$C_m$
such that
$\sup \phi_m^\prime \leq 0$
(this is possible because 
$\phi_m^\prime$
is defined on the compact manifold
$X$
hence bounded from above). Then take the upper semicontinuous upper envelope
$\phi^\prime$
of the
$\phi_m^\prime$
and call
$h$
the (singular) hermitian metric on
$L$
given by the plurisubharmonic weight
$\phi = \phi_L + \phi^\prime$.

\noindent
It is useful to construct the
$\phi_m$
in such a way that 
$\phi^\prime$
has the singularities exactly at
the stable base locus
\[ \mathrm{SBs}(L) := \bigcap_{m \in \mathbb{N}} \mathrm{Bs}(|mL|) \]
of
$L$. This is possible by defining
$\phi_{m+1}$
from
$\phi_m$
as follows: multiply the generators of
$|mm_0L|$
by a section in
$|m_0L|$
and complete this set to a generating set of
$|(m+1)m_0L|$.
By multiplying the completing sections with small positive constants one can
reach around points
$x \in \mathrm{Bs}(|(m+1)m_0L|)$
that
\[ \phi_{m+1} \leq (1-\epsilon_m) \phi_{m} \]
for arbitrarily small
$\epsilon$.
Hence given a positive integer
$M$,
for appropriately chosen
$\epsilon_m$,
there is a constant
$C_x > 0$
such that
$\sup \phi_m \leq C_x \phi_M$.
This implies that 
$\phi^\prime$
has also a singularity in
$x$.

\noindent
The aim is to prove that the Iitaka fibration is (up to birational 
equivalence) the same as the numerically trivial foliation with respect to the
current 
$i\partial\overline{\partial}\phi$.
Since the Iitaka fibration is a fibration this implies in particular
that in this case the numerically trivial foliation is the same as Tsuji's 
numerically trivial fibration.

\noindent
To prove these assertions, first compare Tsuji's and Takayama's intersection 
numbers:
\begin{lem}
With
$L$,
$h$
as above,
\[ (L,h).C \leq ||L,C|| \]
for smooth irreducible curves
$C$
not contained in a Lebesgue zero set.
\end{lem}
\begin{proof}
To begin with, one has to relate the multiplier ideals
$\mathcal{J}(c \cdot |mm_0L|)$
of the linear system
$|mm_0L|$
and the positive rational number 
$c$
with the (analytic) multiplier ideals 
$\mathcal{J}(\phi_m)$.
The ideal 
$\mathcal{J}(c \cdot |mm_0L|)$
is defined via a log resolution, but since
$\phi_m$
is a plurisubharmonic function with analytic singularities defined by 
generating elements of
$|mm_0L|$, 
it follows that
\[ \mathcal{J}(c \cdot mm_0\phi_m) = \mathcal{J}(c \cdot |mm_0L|) \]
by \cite[(5.9)]{Dem00}.

\noindent
As already mentioned in the introduction,
\[ ||L,C|| := \lim_{m \rightarrow \infty} m^{-1} \deg_C mL \otimes 
              \mathcal{J}(||mL||) \]
is defined by using the asymptotic multiplier ideal 
$\mathcal{J}(||mL||)$. 
This ideal is defined to be the unique maximal element among all multiplier 
ideals
$\mathcal{J}(\frac{1}{pm_0} \cdot |pm_0mL|)$
(\cite{Kaw99},\cite{Laz00}). Consequently,
\begin{eqnarray} 
||L,C|| & = & L.C + \lim_{m \rightarrow \infty} m^{-1} \max_{p \in \mathbb{N}} 
                    \deg_C \mathcal{J}(\frac{1}{m_0p}|m_0pmL|) = \\
    & = & L.C + \lim_{m \rightarrow \infty} m^{-1} \lim_{p \rightarrow \infty} 
                    \deg_C \mathcal{J}(m\phi_{mp}) \nonumber \\
    & = & L.C + \lim_{m \rightarrow \infty} m^{-1} \lim_{n \rightarrow \infty} 
                    \deg_C \mathcal{J}(m\phi_{n}) \nonumber . 
\end{eqnarray}

\noindent
The last equality is true because
$ \mathcal{J}(m\phi_{n}) \subset \mathcal{J}(m\phi_{n+1})$
for all
$n$:
The multiplier ideals do not depend on the generating set used to define
$\phi_n$.
By multiplying the generators defining
$\phi_n$
with a section in
$H^0(X, m_0L)$
and completing this set to a generating set of
$H^0(X, (n+1)m_0L)$
it is possible to choose 
$\phi_n \leq \phi_{n+1}$
(as above), hence the inclusion.

\noindent
Next, Tsuji's intersection number may be expressed as
\[ (L,h_n).C = L.C + \limsup_{m \rightarrow \infty} m^{-1}\deg_C 
                     \mathcal{J}(m\phi_n), \]
by \cite[(2.3.1)]{E02} and the fact that
$h_n$
is a metric with analytic singularities, hence restriction to
$C$
and taking the multiplier ideal in the 
$\limsup$
above may be interchanged on smooth curves (\cite[Prop. 2.11]{E02}). An easy
analysis shows that
\begin{eqnarray*}
   \lim_{n \rightarrow \infty} (L,h_n).C & = & L.C + 
   \lim_{n \rightarrow \infty} \limsup_{m \rightarrow \infty} m^{-1}\deg_C 
                     \mathcal{J}(m\phi_n) \\ 
     & \leq & L.C + \lim_{m \rightarrow \infty} m^{-1} 
                  \lim_{n \rightarrow \infty} 
                    \deg_C \mathcal{J}(m\phi_{n}) = ||L;C||. 
\end{eqnarray*} 

\noindent
On the other hand,
$(L,h_n).C = L.C - \sum_{x \in C} \nu(h_{n|C},x)$
by \cite[Prop. 1.2]{E02}. Since the upper semicontinuous upper envelope 
$\phi^\prime$
of the
$\phi_m^\prime$
equals
$\sup_{m} \phi_m^\prime$
outside a set of Lebesgue measure zero (\cite{Lel68}), the envelope of the 
restrictions
$\phi_{m|C}^\prime$
equals almost everywhere the restriction
$(\phi_m^\prime)_{|C}$ 
on all curves outside a Lebesgue zero set. For these curves the lemma follows
from the next statement, using the definition of Lelong numbers via integrals 
(\cite[(2.7)]{Dem00}).
\end{proof}

\begin{lem}
Let
$C \subset X$
be a smooth curve not contained in
$\{ x \in X : \sup_{m} \phi_m^\prime(x) < \phi^\prime(x) \}$.
Then for all
$x \in C$
\[ \lim_{n \rightarrow \infty}\nu(h_{n|C},x) \stackrel{\geq}{\longrightarrow} 
   \nu(h_{|C},x). \]
\end{lem}
\begin{proof}
By definition of Lelong numbers,
$\nu(\phi,x) \geq \nu(\psi,x)$
if
$\phi \leq \psi$. 
Consequently, by the same construction as for the inclusion
$ \mathcal{J}(m\phi_{n}) \subset \mathcal{J}(m\phi_{n+1})$, 
the Lelong numbers 
$\nu(h_{n|C},x)$
of the 
$\phi_n$
form a decreasing sequence of non-negative numbers in every point
$x \in C$
whose limit is 
$\geq \nu(h_{|C},x)$. 
It remains to show the equality:

\noindent
If
$z$
is a local parameter of 
$C$
centered in
$x$,
the restriction
$\phi_n^\prime$
may locally be written as
\[ \phi_n^\prime(z) = \phi_n(z) - \phi_L(z) - C_n =
   \nu(h_{n|C},0) \log|z| + d_n \log(1+\sum_{i=0}^\infty a_i|z|^i) 
   - \phi_L(z) - C_n. \]
For every
$\epsilon > 0$
and a sufficiently small neighborhood of
$0$
it is true that
\[ d_n \log(1+\sum_{i=0}^\infty a_i|z|^i) 
   - \phi_L(z) - C_n \leq -\epsilon\log|z|, \]
hence
$\phi_n^\prime(z) \leq (\nu(h_{n|C},0) - \epsilon) \log|z|$,
which implies
\[ \phi^\prime(z) \leq 
   (\lim_{n \rightarrow \infty} \nu(h_{n|C},0) - \epsilon) \log|z| \]
for almost all
$z$
around 
$0$.
Consequently, 
$\nu(\phi^\prime,0) \geq \lim_{n \rightarrow \infty} \nu(h_{n|C},0) 
 - \epsilon$
for all
$\epsilon > 0$,
and the equality follows.
\end{proof}

\noindent
This already implies that the Iitaka fibration is contained in 
Tsuji's numerically trivial fibration: Take a birational morphism
$\mu: X^\prime \rightarrow X$
from a smooth projective variety 
$X^\prime$
such that the Iitaka fibration induced by the linear system
$|m\mu^\ast L|$
is a morphism 
$f: X \rightarrow Y$
on another smooth variety
$Y$. 
The general fiber of this fibration is smooth. Smooth varieties are numerically
trivial w.r.t. some pair
$(L,h)$
iff 
$(L,h).C = 0$
for all sufficiently general smooth curves in this variety (\cite[3.1]{E02}). 
Hence by the above inequality the numerically trivial fibration w.r.t. 
$(\mu^\ast L, \mu^\ast h)$
contains the Iitaka fibration. By birational equivalence of 
intersection numbers (\cite[2.6]{E02}), the numerically 
trivial fibration w.r.t. 
$(\mu^\ast L, \mu^\ast h)$
is birationally equivalent to that on 
$X$
w.r.t.
$(L,h)$.

\noindent
Next note that there is a positive integer
$m$
such that the Iitaka fibration of
$L$
is induced by the linear system 
$|mL|$
\cite[10.3]{Ii82}.  
\begin{lem}
Let
$L$
be a holomorphic line bundle on a projective complex manifold
$X$
such that
$|mL|$
is a non-empty linear system which induces a rational map
$\phi_{|mL|}: X \dasharrow Y$.
Then
$\phi_{|mL|}$
is the numerically trivial foliation w.r.t. 
$h_{|mL|}$.
\end{lem}
\begin{proof}
By corollary~\ref{numtriv-cor} it is enough to show that for every holomorphic
map
$f: \Delta \rightarrow X$
such that
$\Delta$
is not mapped to a point and does not intersect the base locus of
$|mL|$, 
the unit disk
$\Delta$
is \textbf{not} numerically trivial w.r.t. 
$f^\ast h_{|mL|}$.
But when
$|mL|$
has no  base points in the image of 
$\Delta$,
the metric
$f^\ast h_{|mL|}$
is a smooth metric with smooth positive curvature form different from 
$0$. 
\end{proof}

\begin{lem}
For
$X$,
$L$
and
$h$
as above, let 
$m > 0$
be an integer such that
$|mL|$
is a non-empty linear system,
and 
$f: \Delta \rightarrow X$
a holomorphic map such that
$f^\ast h \neq \infty$,
$f^\ast h_{|mL|} \neq \infty$.
Then
\[ \Delta\ \mathrm{numerically\ trivial\ w.r.t.\ }  f^\ast h
   \Longrightarrow 
   \Delta\ \mathrm{numerically\ trivial\ w.r.t.\ }  f^\ast h_{|mL|}. \]
\end{lem}
\begin{proof}
This is a trivial consequence of the Lelong number inequality
$\nu(f^\ast h,x) \leq \nu(f^\ast h_m,x)$,
see above.
\end{proof}

\noindent
The last lemma implies that the numerically trivial foliation w.r.t. 
$h$
is contained in the numerically trivial foliation w.r.t.
$h_{|mL|}$,
and the lemma before shows that this foliation is the Iitaka fibration
which in turn is contained in Tsuji's numerically trivial fibration by the 
arguments above.

\noindent
\textbf{Remark.} 
This also shows that the Iitaka fibration is the numerically trivial
foliation w.r.t.
$h_{|mL|}$
for an appropriate
$m$.

\section{Bounds for the numerical dimension} \label{BoundND-sec}

\noindent
In this section the ideas of Boucksom and Demailly 
are used to construct a numerically trivial foliation for 
pseudo-effective 
$(1,1)-\!$
classes on compact K\"ahler manifolds. The next paragraph 
tries to collect the scattered and mostly unpublished definitions and 
properties of volumes and moving 
intersection products of pseudo-effective classes without 
claiming any originality or completeness and mostly without proof (in many 
cases they may be 
found in \cite{Bou02}). The main result about the numerically trivial 
foliations will be that the codimension of their leaves 
determines an upper bound for the numerical dimension of the pseudo-effective
class, if the singularities of the foliation are isolated points.

\subsection{Moving intersection numbers of pseudo-effective classes}
\label{MovIntNum-ssec}

\noindent
Starting with Fujita's approximate Zariski decomposition 
(\cite{Fuj94},\cite{DEL00}) Boucksom developped
a notion of volume for arbitrary pseudo-effective classes (\cite{Bou01})
on compact K\"ahler manifolds. This was 
generalized (with small modifications) by Demailly to a ``moving 
intersection product'' of pseudo-effective classes (\cite{Dem02}). This in 
turn allows the definition of a numerical dimension for pseudo-effective 
classes. 

\noindent
Logically one has to start with defining the ``moving intersection numbers'':
\begin{Def}
Let
$X$
be a compact K\"ahler manifold with K\"ahler form
$\omega$.
Let
$\alpha_1, \ldots ,\alpha_p \in H^{1,1}(X,\mathbb{R})$
be pseudo-effective classes and let
$\Theta$
be a closed positive current of bidimension
$(p,p)$.
Then the moving intersection number 
$(\alpha_1 \cdot \ldots \cdot \alpha_p \cdot \Theta)_{\geq 0}$
of the
$\alpha_i$
and
$\Theta$
is defined to be the limit when
$\epsilon > 0$
goes to
$0$
of 
\[ \sup \int_{X-F} (T_1+\epsilon\omega) \wedge \ldots \wedge 
   (T_p+\epsilon\omega) \wedge \Theta \]
where the 
$T_i$'s
run through all currents with analytic singularities in
$\alpha_i[-\epsilon\omega]$,
and 
$F$
is the union of the
$Sing(T_i)$.
\end{Def}

\noindent
It is not difficult to justify the existence of the limit above: First, on
$X-F$ 
the currents
$T_i + \epsilon\omega$
may locally be written as
$T_i + \epsilon\omega = dd^c u_i$
for some bounded plurisubharmonic function
$u_i$. 
By results of Bedford-Taylor \cite{BT76} this implies the existence of the 
integral. In addition Boucksom \cite{Bou01} showed that these integrals are
bounded by a constant only depending on the comological classes
$\{ T_i \}$
and
$\{ \Theta \}$
(this is where the K\"ahler assumption comes in). 
Hence the supremum always exists, and is increasing with increasing
$\epsilon$. 
This implies the existence of the limit. Finally it is easy to 
see that this limit does not depend on the choice of the K\"ahler 
form~$\omega$. 

\noindent
The
$(\alpha_1 \cdot \ldots \cdot \alpha_p \cdot \Theta)_{\geq 0}$
are symmetric in the
$\alpha_i$
and concave and homogeneous in every variable separately. For nef classes
$\alpha_i \in H^{1,1}(X,\mathbb{R})$
the moving intersection number equals the normal cohomological intersection
number
$(\alpha_1 \cdot \ldots \cdot \alpha_p \cdot \{ \Theta \})$
\cite{Bou02}. If some of the pseudo-effective classes coincide one has

\begin{lem} \label{powerint-lem}
For pseudo-effective classes
$\alpha, \alpha_{p+1},\ldots,\alpha_n$
the moving intersection number
$(\alpha^p \cdot \alpha_{p+1} \cdot \ldots \cdot \alpha_n)_{\geq 0}$
is the limit for
$\epsilon \rightarrow 0$
of
\[ \sup \int_{X-F} (T+\epsilon\omega)^p \wedge (T_{p+1}+\epsilon\omega) \wedge
   \ldots \wedge (T_n+\epsilon\omega)\]
where 
$T \in \alpha[-\epsilon\omega]$
and
$T_i \in \alpha_i[-\epsilon\omega]$
have analytic singularities.
\end{lem}
\begin{proof}
See Lemma 3.2.7 in \cite{Bou02}.
\end{proof}

\noindent
\begin{Def}
Let
$X$
be a compact K\"ahler manifold. Then the numerical dimension
$\nu(\alpha)$
of a pseudo-effective class
$\alpha \in H^{1,1}(X,\mathbb{R})$
is defined as
\[ \max \{ k \in \{ 0, \ldots, n\}: (\alpha^k \cdot \omega^{n-k})_{\geq 0} 
                                                                     > 0\}\]
for some (and hence all) K\"ahler classes 
$\omega$.
\end{Def}

\noindent
Now the volume of a pseudo-effective class
$\alpha \in H^{1,1}(X,\mathbb{R})$
on a compact K\"ahler manifold may be defined as a special case of the moving
intersection product:
\[ \mathrm{vol}(\alpha) = (\alpha^n)_{\geq 0}. \]
But there are other useful possibilities to define it: First remember that
Fujita considered projective 
$n-\!$
dimensional algebraic varieties
$X$
and line bundles
$L$
over 
$X$,
and defined
the volume of 
$L$
by
\[ \mathrm{vol}(L) := 
                 \limsup_{k \rightarrow +\infty} \frac{n!}{k^n} h^0(X,kL). \]
If
$L$
is nef
the volume of 
$L$
is the self-intersection
$L^n$,
by Riemann-Roch and 
$h^q(X,kL) \sim O(k^{n-q})$
(\cite[(6.7)]{Dem00}). For arbitrary pseudo-effective classes
$\alpha \in H^{1,1}(X, \mathbb{R})$
on compact K\"ahler manifolds
$X$
Boucksom generalized this volume by defining
\[ \mathrm{vol}(\alpha) = \sup \int_X T_{ac}^n \]
where the supremum is taken over all closed positive 
$(1,1)-\!$
currents
$T$
with
$\{ T \} = \alpha$
and
$T_{ac}$
is the absolute continuous part of the Lebesgue decomposition
$T = T_{ac} + T_{sg}$.
Again, the K\"ahler assumption is necessary to guarantee that
$T_{ac}^n$
is locally integrable. By using singular Morse inequalities and the 
Calabi-Yau theorem Boucksom proved that
$\mathrm{vol}(L) = \mathrm{vol}(c_1(L))$
and that
$\mathrm{vol}(L) > 0$
iff 
$L$
is a big line bundle, i.e. iff there is a closed strictly positive current 
representing
$c_1(L)$.

\noindent
Note that it is not necessary to look at all closed positive 
$(1,1)-\!$
currents for taking the supremum. This is a consequence of an approximation 
theorem of Demailly: 
\begin{thm}[\cite{Dem92}] \label{DemApp-thm}
Let
$T = \theta + dd^c \phi$
be a closed almost positive 
$(1,1)-\!$
current on a complex manifold 
$X$
with hermitian metric
$\omega$
such that 
$\theta$
is a smooth form. Suppose that 
$T \geq \gamma$
for some real 
$\mathcal{C}^\infty-\!$
form
$\gamma$. 
Then there exists a decreasing sequence 
$\phi_k$
of almost plurisubharmonic functions with analytic singularities such that the
$T_k := \theta + dd^c \phi_k$
verify
\begin{itemize}
\item[(i)]
The 
$\phi_k$ 
converge pointwise and
$L^1_{loc}$
against
$\phi$,
hence the 
$T_k$ 
converge weakly against
$T$.
\item[(ii)]
$T_k \geq \gamma - \epsilon_k\omega$
for some sequence of positive numbers
$\epsilon_k \rightarrow 0$.
\item[(iii)]
The Lelong numbers
$\nu(T_k,x)$
converge uniformly against
$\nu(T,x)$
w.r.t. 
$x \in X$.
\end{itemize}
\end{thm} 

\noindent
Using another approximation theorem (\cite{Dem82}) Boucksom 
slightly modified this statement (\cite{Bou01}): 
\begin{thm}
Let the assumptions and notations be the same as in the theorem before. Then
there exists a decreasing sequence 
$\phi_k$
of plurisubharmonic functions with analytic singularities such that the
$T_k := \theta + dd^c \phi_k$
verify
\begin{itemize}
\item[(i)]
The 
$T_k$ 
converge weakly against
$T$,
and
$T_{k,ac} \rightarrow T_{ac}$
almost everywhere.
\item[(ii)]
$T_k \geq \gamma - \epsilon_k\omega$
for some sequence of positive numbers
$\epsilon_k \rightarrow 0$.
\item[(iii)]
The Lelong numbers
$\nu(T_k,x)$
converge uniformly against
$\nu(T,x)$
w.r.t. 
$x \in X$.
\end{itemize}
\end{thm}

\noindent
So one may define instead
\[ \mathrm{vol}(\alpha) = 
                   \lim_{\epsilon \rightarrow 0^+} \sup \int_X T_{ac}^n \]
where the
$T$'s
run through all closed 
$(1,1)-\!$
currents with analytic singularities in
$\alpha[-\epsilon\omega]$,
that is
$\{ T \} = \alpha$
and
$T \geq -\epsilon\omega$
for some hermitian metric
$\omega$
on
$X$.

\noindent
Here, closed
$(1,1)-\!$
currents with analytic singularities are currents whose almost 
plurisubharmonic potentials locally look like
\[ \frac{\alpha}{2}\log(|f_1|^2 + \ldots + |f_p|^2) \]
with
$f_1, \ldots, f_n$
holomorphic, up to a bounded
$\mathcal{C}^\infty-\!$
function. Such currents 
$T$
are particularly useful because their absolut 
continuous part is the same as the residual part
$R$
in the Siu-decomposition
$T = \sum_i a_i[D_i] + R$. 
Consequently, one may compute
$\int_X T_{ac}^n$
by blowing up the (integral closure) of the ideal of singularities locally
generated by the
$f_i$
and integrating the smooth form given by the pull back of
$T$
minus the integration currents of the exceptional divisors as they occur in 
the
inverse image of the singularity ideal. In Fujita's setting this corresponds
to blowing up the base locus of the multiples
$mL$
and decomposing the pull back of
$L$
into an effective part
$E_m$
and a free part  
$D_m$, and Fujita's theorem \cite[(14.6)]{Dem00} tells us that
\[ \mathrm{vol}(L) = \lim_{m \rightarrow \infty} D_m^n. \]

\noindent
Finally, the last definition of
$\mathrm{vol}(\alpha)$
is equivalent to the first one, with moving intersection numbers, by
lemma~\ref{powerint-lem}.

\subsection{Numerical triviality for pseudo-effective classes}

\noindent
First repeat and codify the informal definitions of numerical triviality and
numerically trivial foliations w.r.t. a pseudo-effective class from the 
introduction:
\begin{Def} \label{NumTrivPseff-def}
Let
$X$
be a compact K\"ahler manifold with K\"ahler form 
$\omega$
and pseudo-effective class 
$\alpha \in H^{1,1}(X,\mathbb{R})$.
A submanifold
$Y \subset X$
(closed or not) is \textbf{numerically trivial} w.r.t. 
$\alpha$
iff for every immersed disk 
$\Delta \subset Y$,
\[ \lim_{\epsilon \downarrow 0} \sup_T \int_{\Delta^\prime - \mathrm{Sing}\ T} 
    (T + \epsilon \omega)  = 0, \]
where the
$T$'s
run through all currents with analytic singularities in
$\alpha[-\epsilon \omega]$
and
$\Delta^\prime = \{ t: |t| < 1 - \delta \}$
is any smaller disk contained in 
$\Delta = \{ t: |t| < 1 \}$. 
\end{Def} 

\noindent
As a convention set
$\int_{\Delta - \mathrm{Sing}\ T} 
    (T + \epsilon \omega)  = 0$
if
$\Delta - \mathrm{Sing}\ T = \emptyset$.
Furthermore note that the restriction to disks
$\Delta^\prime$
may be replaced by the assumption that it is possible to continue the 
immersion 
$\Delta \subset Y$
holomorphically.

\begin{Def}
Let
$X$
be a compact K\"ahler manifold with a pseudo-effective class
$\alpha \in H^{1,1}(X, \mathbb{R})$. 
A foliation
$\{ \mathcal{F},(U_i,p_i) \}$
is \textbf{numerically trivial} w.r.t.
$\alpha$
iff 
\begin{itemize}
\item[(i)]
every fiber of
$p_i$
is numerically trivial w.r.t. 
$\alpha$,
\item[(ii)]
and if 
$\Delta^2 \hookrightarrow U_i$
is an immersion such that the projection onto the first coordinate coincides 
with the projection
$p_i: U_i \rightarrow \Delta^{n-k}$,
then for any
$\Delta^\prime \subset\subset \Delta$
and any sequence of currents
$T_k \in \alpha[-\epsilon_k\omega]$,
$\epsilon_k \rightarrow 0$, 
the integrals
$\int_{(\{z_1 = a\} \cap \Delta^\prime) - \mathrm{Sing\ }T_k} 
  (T_k + \epsilon_k\omega)$
are uniformly (in
$a$)
bounded from above.
\end{itemize}
\end{Def}

\noindent
Note that no exceptional fibers are allowed: if the fibers are completely 
contained in the common singularity locus of the
$T \in \alpha[-\epsilon\omega]$,
then they are numerically trivial by the convention above, otherwise the limit
in definition~\ref{NumTrivPseff-def} is supposed to be
$0$.
The uniform boundedness is essential for the proof of the Local Key Lemma 
below.

\noindent
To construct a maximal numerical trivial foliation w.r.t. this notion, it is enough to prove an analog for the 
Local~Key~Lemma~\ref{LKlem}:
\begin{lem}[\textbf{Local Key Lemma for pseudo-effective classes}]
Let
$X$
be a compact K\"ahler manifold with a pseudo-effective class
$\alpha \in H^{1,1}(X, \mathbb{R})$. 
Let
$W \cong \Delta^n$
be an open subset of
$X$
with a projection
$p: W \rightarrow \Delta^k$
onto the last 
$k$
factors, and let
$V = \{ z_1 = \ldots = z_{n-k} = 0 \}$
be a complex submanifold of
$W$.  
If every fiber of
$p$
and also
$Y$
are numerically trivial w.r.t. 
$\alpha$,
then 
$W$
will also be numerically trivial w.r.t. 
$\alpha$.
\end{lem}

\noindent
Then the maximal  numerically trivial foliation w.r.t. 
$\alpha$
may be constructed in the same way as in section~\ref{MaxNTFol-ssec}.

\noindent
The proof of this Local Key Lemma for pseudo-effective classes imitates the 
proof of the Local Key Lemma for closed positive
$(1,1)-\!$
currents: There, the numerical triviality of the 
fibers of the projection implies that the residue current of the Siu
decomposition is a pull 
back of a current on the base (see the Pullback~Lemma~\ref{pullback-lem}). Of 
course, the Pullback~Lemma is not true for 
pseudo-effective classes. But it is enough to prove that
the restriction onto different horizontal sections are quite the same, hence 
the numerical triviality of
$V$
implies the numerical triviality of all horizontal sections, hence that of
$W$. 
This argument is made exact by
\begin{prop} \label{BdLKL-prop}
Let
$X$
be a compact K\"ahler manifold with K\"ahler form
$\omega$,
and let
$T_k = T_k^\prime + \epsilon_k\omega$,
$\epsilon_k \rightarrow 0$, 
be a sequence of closed positive
$(1,1)$-\!
currents on
$X$
such that the 
$T_k^\prime$
represent the same cohomology class.
Let
$\Delta^2 \hookrightarrow X$
be an immersion (with coordinates
$z_1, z_2$).
Let
$\Delta^\prime \subset\subset \Delta$
be a disk, and consider the functions
$f_k: \Delta^\prime \rightarrow \mathbb{R}^+, a \mapsto 
      \int_{(\{z_1=a\}\cap\Delta^\prime)-\mathrm{Sing\ } T_k} T_k$
and
$g_k: \Delta^\prime \rightarrow \mathbb{R}^+, b \mapsto 
      \int_{(\{z_1=b\}\cap\Delta^\prime)-\mathrm{Sing\ } T_k} T_k$.
Suppose that 
$\lim_{k \rightarrow \infty} f_k(a) = 0$
for all
$a \in \Delta$,
and that the 
$f_k$
are uniformly (in
$a$)
bounded from above. Suppose furthermore that
$\lim_{k \rightarrow \infty} g_k(0) = 0$.
Then
$\lim_{k \rightarrow \infty} g_k(b) = 0$
for all
$b \in \Delta^\prime$,
and the 
$g_k$
are uniformly (in
$b$)
bounded from above.
\end{prop}
\begin{proof}
Since the integrals are always evaluated outside the singularities of
$T_k$,
and since the mass of the integration current of a divisor is always 
concentrated in the divisor, one can assume without loss of generality that the
Siu decomposition of
$T_k$
does not contain any integration currents of divisors. Consequently, 
$T_k$
has only finitely many isolated singularities on any compact subset of 
$p_1^{-1}(\Delta^\prime)$
where
$\Delta^\prime \subset\subset \Delta$
is any disk and
$p_1: \Delta^2 \rightarrow \Delta$
is the projection onto the first coordinate, and
$T_k$
may be written on 
$p_1^{-1}(\Delta^\prime)$
as
\[ T_k = \theta_{11}^k idz_1 \wedge d\overline{z}_1 +
         \theta_{12}^k idz_1 \wedge d\overline{z}_2 + 
         \theta_{21}^k idz_2 \wedge d\overline{z}_1 +
         \theta_{22}^k idz_2 \wedge d\overline{z}_2,  \]
where the 
$\theta_{ij}^k$
are smooth functions outside these singularities, and integrable on
$\Delta^2$.
That 
$T_k$
is a real current implies
$\theta_{ij}^k  = \overline{\theta_{ji}^k}$.

\noindent
To prove the proposition it is enough to show that
\[ \lim_{k \rightarrow \infty} |\int_{\Delta^\prime_b - \mathrm{Sing\ }T_k}T_k 
       - \int_{\Delta^\prime_0 - \mathrm{Sing\ }T_k} T_k| = 0. \]
for a sequence of disks
$\Delta^\prime \subset\subset \Delta$
exhausting 
$\Delta$
(where
$\Delta^\prime_b = \{z_2 = b\} \cap \Delta^\prime$).
Now, choose a path
$\gamma \in \Delta$
from
$0$
to
$b$.
Then,
\[ \begin{array}{rcl}
   |\int_{\Delta^\prime_b - \mathrm{Sing\ }T_k} T_k - 
    \int_{\Delta^\prime_0 - \mathrm{Sing\ }T_k} T_k| & = &
   |\int_{\Delta^\prime} (\theta_{11}^k(z_1,b) - \theta_{11}^k(z_1,0)) 
                idz_1 \wedge d\overline{z}_1| 
    \end{array} \]
equals (by Stokes and Fubini)
\[ \begin{array}{rcl}
   |\int_{\Delta^\prime} (\int_\gamma d\theta_{11}^k) 
          idz_1 \wedge d\overline{z}_1|
    & = & |\int_{\Delta^\prime \times \gamma} d(\theta_{11}^k 
               idz_1 \wedge d\overline{z}_1)|. 
   \end{array} \]
Since the closedness of
$T$
implies
\[ d(\theta_{11}^k idz_1 \wedge d\overline{z}_1) =
      - d(\theta_{12}^k idz_1 \wedge d\overline{z}_2 + 
          \theta_{21}^k idz_2 \wedge d\overline{z}_1 +
          \theta_{22}^k idz_2 \wedge d\overline{z}_2),\]
this integral equals by Stokes
\[ |\int_{\partial(\Delta^\prime \times \gamma)} 
        (\theta_{12}^k idz_1 \wedge d\overline{z}_2 + 
          \theta_{21}^k idz_2 \wedge d\overline{z}_1 +
          \theta_{22}^k idz_2 \wedge d\overline{z}_2)|, \]
and since
$z_2$
is constant on
$\Delta^\prime \times \partial \gamma$,
this simplifies to
\[   |\int_{(\partial\Delta^\prime) \times \gamma} 
        (\theta_{12}^k idz_1 \wedge d\overline{z}_2 + 
          \theta_{21}^k idz_2 \wedge d\overline{z}_1 +
          \theta_{22}^k idz_2 \wedge d\overline{z}_2)|.   \]

\noindent
Observe that these integrals do not depend on the chosen path
$\gamma$.
Consequently, cover the disk
$\Delta_{0,b}$
with center in
$b/2$
and radius
$|b/2|$
with a family of paths
$\gamma_a$
from
$0$
to
$b$.
Then to prove
$\lim_{k \rightarrow \infty} 
 |\int_{\Delta_b^\prime} T_k - \int_{\Delta_0^\prime} T_k| = 0$
it is enough to show that
\[ \lim_{k \rightarrow \infty} 
  \int_a |\int_{(\partial\Delta^\prime) \times \gamma_a} 
        (\theta_{12}^k idz_1 \wedge d\overline{z}_2 + 
          \theta_{21}^k idz_2 \wedge d\overline{z}_1 +
          \theta_{22}^k idz_2 \wedge d\overline{z}_2))|da = 0. \]
The term with
$\theta_{22}^k$
vanishes since 
$idz_2 \wedge d\overline{z}_2$
is pulled back to
$0$
in any chart of
$(\partial\Delta^\prime) \times \gamma_a$. 
Since
$\theta_{12}^k = \overline{\theta_{21}^k}$
the remaining integral may be bounded from above by
\[ C \cdot \int_{\partial\Delta^\prime \times \Delta_{0,b}}
            |\theta_{12}^k|dV, \]
where
$C$
is independent of
$b$
and
$k$,
and
$dV$
is a volume element on
$\partial\Delta^\prime \times \Delta_{0,b}$.

\noindent
Now interpret
$T_k$
as a semipositive hermitian form 
$\langle .\ , . \rangle$
on every tangent space
$T_{X,x}$
(where
$T$
has no singularities). Then the Schwarz inequality implies that
\[ |\theta_{12}^k| = 
   |\langle \frac{\partial}{\partial \overline{z}_1},
           \frac{\partial}{\partial \overline{z}_2} \rangle| \leq
   |\langle \frac{\partial}{\partial \overline{z}_1},
            \frac{\partial}{\partial \overline{z}_1} \rangle|^{\frac{1}{2}} 
   \cdot 
   |\langle \frac{\partial}{\partial \overline{z}_2},
            \frac{\partial}{\partial \overline{z}_2} \rangle|^{\frac{1}{2}} = 
   |\theta_{11}^k|^{\frac{1}{2}} \cdot |\theta_{22}^k|^{\frac{1}{2}}. \]
Hence the integral above is
$\leq$
the square root of the product
\[ \int_{\partial\Delta^\prime \times \Delta_{0,b}}
            |\theta_{11}^k| dV \cdot
   \int_{\partial\Delta^\prime \times \Delta_{0,b}}
            |\theta_{22}^k| dV, \]
again by the Schwarz inequality.

\noindent
\textit{Claim.}
There exists a bound
$M^\prime > 0$
such that for all
$k$
there is a disk
$\Delta_k^\prime \subset\subset \Delta$
containing
$\Delta^\prime$
with
\[ \int_{\partial\Delta_k^\prime \times \Delta^\prime}
            |\theta_{11}^k| dV <
   M^\prime. \]
\begin{proof}
Suppose that
$\Delta^\prime \subset\subset \Delta^{\prime\prime} \subset\subset \Delta$,
and look at the
$(1,1)$-\!
form
$\eta = idz_2 \wedge d\overline{z}_2$.
There exists a
$C > 0$,
such that
$\eta \leq C \cdot \omega$
on
$\Delta^{\prime\prime} \times \Delta^\prime$.
Hence,
\[ \int_{(\Delta^{\prime\prime}-\Delta^\prime) \times \Delta^\prime} 
            |\theta_{11}^k| dV =
   \int_{(\Delta^{\prime\prime}-\Delta^\prime) \times \Delta^\prime} 
   (T_k^\prime + \epsilon_k\omega) \wedge \eta \leq
   C \cdot \int_X (T_k^\prime + \epsilon_k\omega) \wedge \omega, \]
and the last integral only depends on the cohomology class of
$T_k^\prime$
(and
$\omega$).
By Fubini one gets a disk
$\Delta_k^\prime$
as above.
\end{proof}

\noindent
For the second term note that the assumptions on the functions 
$f_k$
imply
$\lim_{k \rightarrow \infty} \int_{\Delta^\prime} f_k ida \wedge d\overline{a} 
 = 0$,
by Lebesgue's dominated convergence, 
and the measure of the sets
$\{a: f_k(a) > \delta\}$
tends to
$0$,
too, for
$k \rightarrow \infty$.

\noindent
Hence, as above, for a given
$\epsilon > 0$
it is possible to bound the measure of
$\{a: f_k(a) > \delta\}$
small enough such that for all
$k$
big enough there is a disk
$\Delta_k^{\prime\prime} \subset\subset \Delta$
containing 
$\Delta^\prime$
with
\[ \int_{\partial\Delta_k^{\prime\prime} \times \Delta^\prime}
            |\theta_{22}^k| dV <
   \epsilon. \]
Choosing
$\delta$
small enough and
$M^\prime$
big enough (but both independent of
$k$!)
one can assume that the two disks
$\Delta_k^\prime$
and
$\Delta_k^{\prime\prime}$
coincide (at least for
$k$
big enough). Since
$M^\prime$ 
is independent of
$\epsilon$, 
the difference
$\int_{\Delta_{k,b}} T_k - \int_{\Delta_{k,0}} T_k$
tends to
$0$
for
$k \rightarrow \infty$,
and uniformly in
$b$.
Since
$\int_{\Delta_{k,0}} T_k \stackrel{k \rightarrow \infty}{\longrightarrow} 0$,
this is also true for
$\int_{\Delta_b^\prime} T_k - \int_{\Delta_0^\prime} T_k$.
Consequently,
$\lim_{k \rightarrow \infty} g_k(a) = 0$,
and the uniformity in
$b$
implies the uniform boundedness of the
$g_k$.
\end{proof}

\noindent
\textit{Proof of the Local Key Lemma for pseudo-effective classes.}  If
$\Delta$
is a disk immersed in
$W$
such that 
$p$ 
projects it on a point in
$\Delta^k$, 
there is nothing to prove. 

\noindent
If
$\Delta$
is a disk immersed in
$W$
not intersecting
$Y$
which is projected biholomorphically onto
$\Delta^k$,
then a coordinate change and further cutting down leads
to the configuration described in the proposition.
Note that it is sufficient to check on any disk
$\Delta^{\prime} \subset\subset \Delta$
that 
\[ \lim_{k \rightarrow \infty}  
   \int_{\Delta^{\prime} - \mathrm{Sing} T_k} T_k + \frac{1}{k} \omega = 0 \]
for arbitrary 
sequences 
$T_k$
of currents with analytic singularities in
$\alpha[- \frac{1}{k}\omega]$. 
The assumptions of the Local Key Lemma imply that
\[ \lim_{k \rightarrow \infty} \int_{\{z_1 = a\} - \mathrm{Sing} T_k} 
   T_k + \frac{1}{k}\omega = \lim_{k \rightarrow \infty} f_k(a) = 0 \]
for all
$a$
and
$\lim_{k \rightarrow \infty} \int_{\{z_2 = 0\} - \mathrm{Sing} T_k} 
 T_k + \frac{1}{k}\omega = 0$.
The definition of a numerically trivial foliation implies the uniform 
boundedness of the
$f_k$,
so it is possible to apply the proposition.

\noindent
If
$\Delta$
is a disk immersed in
$W$
not satisfying one of the two conditions above, then for any
$\Delta^\prime \subset\subset \Delta$
there are disks
$\Delta_i^{\prime\prime} \subset\subset \Delta_i^\prime \subset \Delta$
such that
$\bigcup \Delta_i^{\prime\prime} \supset \Delta^\prime$
(hence it is enough to consider finitely many of these disks), and
there are projections
$p_i: W \rightarrow \Delta^{n-k}$
(possibly different from
$p$)
such that the restriction onto
$\Delta_i^\prime$
is a submersion. Since the fibers and sections of these 
$p_i$
are composed of disks already shown to be numerically trivial, it is possible 
to apply again the proposition on
$\Delta_i^{\prime\prime} \subset\subset \Delta_i^\prime$
(by possibly further cutting down and a coordinate change).
Since there are only finitely many
$i$'s,
$\Delta^\prime$
is also numerically trivial.

\noindent
Finally, the uniform boundedness property of the foliation follows directly 
from the uniform boundedness shown in the proposition.
\hfill $\Box$

\subsection{The Iitaka fibration and the nef fibration}

\label{ITnefFib-ssec} 

\noindent
A remarkable fact about the construction of a numerically trivial foliation
w.r.t. a pseudo-effective class
$\alpha$
is that it works also if one restricts to non-empty subsets of currents with
analytic singularities in
$\alpha[-\epsilon\omega]$.
\begin{lem}
Let
$X$
be a compact K\"ahler manifold with K\"ahler form
$\omega$
and
$\Theta$
a closed positive
$(1,1)-\!$
current representing the cohomology class
$\alpha \in H^{1,1}(X,\mathbb{R})$.
Then the foliation constructed w.r.t. the subsets 
$\{\Theta\} \subset \alpha[-\epsilon\omega]$
is the numerically trivial foliation
w.r.t.
$\Theta$.
\end{lem}
\begin{proof}
Comparing definitions~\ref{NumTrivSubVar-def} and \ref{NumTrivPseff-def}, and 
taking into account the 
criterion~\ref{NumTrivCrit-prop} for numerical triviality one immediately gets
that the numerically trivial foliation
w.r.t.
$\Theta$
is contained in that constructed w.r.t. the subsets 
$\{\Theta\} \subset \alpha[-\epsilon\omega]$.
But the other inclusion is also not difficult to prove: Every holomorphic map
$\Delta \rightarrow X$
maps 
$\Delta$
onto a 
$1-\!$
dimensional analytic subset, and the integrals in 
definition~\ref{NumTrivPseff-def} may be 
taken outside the singularities of this set. 
\end{proof}

\noindent
It is also clear that this foliation contains the numerically trivial 
foliation w.r.t.
$\alpha$. 
In particular:
\begin{prop}
Let
$X$
be a projective manifold and
$L$
a nef line bundle on
$X$
such that the Kodaira-Iitaka dimension
$\kappa(L) \geq 0$.
Then the nef foliation of
$L$
is contained in the Iitaka fibration. 
\hfill $\Box$
\end{prop} 

\noindent
In analogy to Tsuji's numerically trivial fibration one can define the
pseudo-effective fibration of a pseudo-effective line bundle 
$L$
as the maximal fibration contained in the numerically trivial foliation w.r.t.
$c_1(L)$.
\begin{prop}
Let
$X$
be a projective manifold and
$L$
a nef line bundle on
$X$.
Then the nef fibration of~\cite{BCEKPRSW00} is equal to the pseudo-effective 
fibration.
\end{prop}
\begin{proof}
This is just a consequence of the definitions: A curve 
$C$
which is numerically trivial w.r.t. 
$c_1(L)$
satisfies
$(L.C)_{\geq 0} = L.C = 0$,
and vice versa.
\end{proof}

\noindent
To summarize, all this gives a (sufficient) geometric reason that the fibers 
of the nef 
fibration are strictly contained in the fibers of the Iitaka 
fibration: this happens if the nef foliation is not a fibration. It would be 
interesting to decide if the converse is also true.

\subsection{Currents with minimal singularities} \label{MinCurr-ssec}

\noindent
To state and to prove the results about upper bounds for the numerical 
dimension of
a pseudo-effective class, a further notion is still missing: that of currents 
with minimal singularities.
\begin{Def}
Let
$\phi_1$
and
$\phi_2$
be two almost plurisubharmonic functions on a complex manifold
$X$.
Then 
$\phi_1$
is said to be less singular than
$\phi_2$
in 
$x \in X$
iff
\[ \phi_2 \leq \phi_1 + O(1) \]
in a neighborhood of 
$X$.
The fact that
$\phi_1$
is less singular than
$\phi_2$
in every point is denoted by
$\phi_1 \preceq \phi_2$. 
\end{Def}

\noindent
Now let
$X$
be compact K\"ahler
and
$\alpha \in H^{1,1}(X,\mathbb{R})$.
Let
$\theta$
be a smooth
$(1,1)-\!$
form representing
$\alpha$.
Then every current in
$\alpha$
may be written as
$T = \theta + dd^c\phi$
for some almost plurisubharmonic function
$\phi$ 
and
\[ T_1 \preceq T_2 \]
shall denote the fact that
$\phi_1 \preceq \phi_2$.
\begin{prop}
Let
$\gamma$
be a smooth 
$(1,1)-\!$
form on
$X$. 
Every non-empty subset of 
$\alpha[\gamma]$
admits a lower bound in
$\alpha[\gamma]$ 
w.r.t. 
$\preceq$.
\end{prop}
\begin{proof}
The proof is almost trivial and of course contained in \cite{DPS00} but is 
repeated for emphasizing a certain uniqueness property.

\noindent
Let
$(T_i)_{i \in I}$
be the given subset of
$\alpha[\gamma]$.
Write 
$T_i = \theta + dd^c\phi_i$
where
$\phi_i$
is almost plurisubharmonic and
$dd^c\phi_i \geq \gamma - \theta$.
Since 
$X$
is compact, all almost plurisubharmonic functions are bounded from above hence 
one may suppose that
$\phi_i \leq 0$
by subtracting a constant. If one choose this constant such that 
$\sup_{x \in X} \phi_i(x) = 0$
the
$\phi_i$
will be unique: An almost plurisubharmonic function
$\phi$
with
$dd^c \phi = 0$
is a holomorphic function.

\noindent
The 
$\phi_i$
have an almost plurisubharmonic upper envelope 
$\phi$
such that
$\theta + dd^c \phi \in \alpha[\gamma]$.
The current
$T = \theta + dd^c \phi$
is obviously a lower bound for the
$(T_i)_{i \in I}$,
with the following property: If 
$S \preceq T_i$
for all
$I$,
then
$S \preceq T$.  
\end{proof}

\noindent
\textbf{Remark.}
The construction above shows that this lower bound
$T = T_{min}$
is unique only up to
$L^\infty$.
On the other hand, given the smooth
$(1,1)-\!$
form
$\theta$
in 
$\alpha$, 
the construction leads to a well defined current 
$T_{min} = \theta + dd^c \phi_{min}$
via the upper envelope. Here, the almost plurisubharmonic function 
$\phi_{min}$ 
satisfies
$\phi_i \leq \phi_{min}$
where the 
$\phi_i$
are chosen as above.

\noindent
This current will be used in the following.

\smallskip

\noindent
The currents with minimal singularities may be used to define minimal 
multiplicities of pseudo-effective classes, having a look at Boucksom's 
construction of higher dimensional Zariski 
decompositions~\cite{Bou02b}. In this paper, he interpreted the Lelong numbers 
of a current 
$T_{min,\epsilon}$
with minimal singularities in
$\alpha[-\epsilon\omega]$
as the obstructions to reach smooth currents in
$\alpha[-\epsilon\omega]$.
This led him to
\begin{Def}
The minimal multiplicity of a pseudo-effective class 
$\alpha \in H^{1,1}(X,\mathbb{R})$
in
$x \in X$
is defined as
\[ \nu(\alpha,x) := \sup_{\epsilon > 0}\nu(T_{min,\epsilon},x).\]
The generic minimal multiplicity on a prime divisor
$D \subset X$
is defined as
\[ \nu(\alpha,D) := \inf_{x \in D}\nu(\alpha,x).\]
\end{Def}

\noindent
Denoting by 
$T_{min}$
a current with minimal singularities in
$\alpha[0]$
one has always
\[ \nu(\alpha,x) \leq \nu(T_{min},x),\ \nu(\alpha,D) \leq \nu(T_{min},D). \]
There are examples where
$\nu(\alpha,D) < \nu(T_{min},D)$, see section~\ref{firstex-ssec}.

\smallskip

\noindent
The following approximation of 
$T_{\mathrm{min}}$
which will be useful later on:
\begin{thm} \label{NumdimApp-thm}
Let
$X$ 
be a compact K\"ahler manifold with K\"ahler class
$\omega$,
let
$\alpha \in H^{1,1}(X,\mathbb{R})$
be a pseudo-effective class. 
Then there exists a sequence of closed
$(1,1)-\!$
currents
$T_k$
with analytic singularities in
$\alpha[-\epsilon_k\omega]$
for some sequence
$(\epsilon_k) \rightarrow 0$
of positive real numbers such that
\begin{enumerate}
\item[(i)]
the 
$T_k$
converge weakly against a closed positive
$(1,1)-\!$
current
$T$
which has minimal singularities
in
$\alpha[0]$,
\item[(ii)]
$\nu(T_k,x) \rightarrow \nu(\alpha,x)$
for every point
$x \in X$,
\item[(iii)]
for all
$i$
\[ \int_{X-Sing(T_k)} (T_k+\epsilon_k\omega)^p \wedge \omega^{n-p} 
   \rightarrow 
   (\alpha^{p}.\omega^{n-p})_{\geq 0} . \]
\end{enumerate}
\end{thm}
\begin{proof}
To compute
$(\alpha^{p}.\omega^{n-p})_{\geq 0}$
it is enough to determine the limit of the
\[ s_\epsilon := \sup_T \int_{X-Sing(T)} (T+\epsilon\omega)^p \wedge 
   \omega^{n-p} \]
where
$T \in \alpha[-\epsilon\omega]$
has analytic singularities, by Lemma~\ref{powerint-lem}. Consequently, for 
each
$p$
there is a sequence of closed 
$(1,1)-\!$
currents 
$(T_k^{(p)})_{k \in \mathbb{N}}$
with analytic singularities such that
$T_k^{(p)} \in  \alpha[-\epsilon_k\omega]$
for some sequence
$\epsilon_k \rightarrow 0$
of positive real numbers and
\[ \int_{X-Sing(T_k^{(p)})} (T_k^{(i)}+\epsilon_k\omega)^p \wedge 
    \omega^{n-p}
    \rightarrow (\alpha^{p}.\omega^{n-p})_{\geq 0} .\] 

\noindent
Now let
$\theta$
be a smooth
$(1,1)-\!$
form on 
$X$
representing
$\alpha$.
Let
$T_{min,k} = \theta + dd^c\phi_{min,k}$
be the current with minimal singularities in
$\alpha[-\epsilon_k\omega]$
associated to
$\theta$,
as described in the remark above. Since
$T_k^{(p)} = \theta + dd^c\phi_k^{(p)} \in \alpha[-\epsilon_k\omega]$
this implies
$\phi_k^{(p)} \leq \phi_{min,k} \leq 0$.
Furthermore the
$T_{min,k}$
converge weakly against a current
$T_{min}$
with minimal singularities in
$\alpha[0]$.

\noindent 
By Demailly's Approximation Theorem~\ref{DemApp-thm} there exists a decreasing
sequence of almost plurisubharmonic functions
$\phi_{k,l}$
with analytic singularities converging pointwise and
$L^1_{loc}$
against
$\phi_{min,k}$
such that
$T_{k,l} = \theta + dd^c \phi_{k,l} \in \alpha[-\epsilon_{k,l}\omega]$
for some sequence
$(\epsilon_{k,l})_{l \in \mathbb{N}} \stackrel{>}{\rightarrow} \epsilon_k$
of positive real numbers. Furthermore
$\nu(T_{k,l},x) \stackrel{<}{\rightarrow} \nu(T_{min,k},x)$
for every point
$x \in X$.

\noindent
Let
$\mu : Y \rightarrow X$
be a common resolution of the singularities of
$T_{k,l}$
and the
$T_k^{(p)}$. 
Then 
\[ \mu^\ast T_k^{(p)} = R_k^{(p)} + [D_k^{(p)}], 
   \mu^\ast T_{k,l} = R_{k,l} + [D_{k,l}] \]
where
$R_k^{(p)}, R_{k,l}$
are smooth and
$D_k^{(p)}, D_{k,l}$
are effective 
$\mathbb{R}-\!$
divisors. Since the 
$\phi_{k,l}$
form a decreasing sequence,
$\phi_k^{(p)} \leq \phi_{k,l}$
and 
$T_{k,l}$
is less singular than
$T_k^{(p)}$.
In particular
$D_{k,l} \leq D_k^{(p)}$,
hence the \textbf{class}
$\{R_{k,l} - R_k^{(i)}\} = \{D_k^{(i)} - D_{k,l}\}$
is pseudo-effective. Consequently,
\[ \int_Y (R_{k,l}+\epsilon_{k,l}\mu^\ast\omega) \wedge 
          (R_k^{(p)}+\epsilon_{k,l}\mu^\ast\omega)^{p-1} \wedge 
          \mu^\ast\omega^{n-p} \geq
   \int_Y (R_k^{(p)}+\epsilon_{k,l}\mu^\ast\omega)^p \wedge 
          \mu^\ast\omega^{n-p}, \]
since the integrals over the compact manifold
$Y$
only depend on the cohomology classes, and all factors besides
$R_{k,l}+\epsilon_{k,l}\mu^\ast\omega$
and
$R_k^{(p)}+\epsilon_{k,l}\mu^\ast\omega$
are smooth. Iterating gives 
\[ \int_Y (R_{k,l}+\epsilon_{k,l}\mu^\ast\omega)^p \wedge 
          \mu^\ast\omega^{n-p} \geq
   \int_Y (R_k^{(p)}+\epsilon_{k,l}\mu^\ast\omega)^p \wedge 
          \mu^\ast\omega^{n-p}. \]
Noting that
\[ \int_Y (R_{k,l}+\epsilon_{k,l}\mu^\ast\omega)^p \wedge 
          \mu^\ast\omega^{n-p} = 
   \int_{X-Sing(T_{k,l})} (T_{k,l} + \epsilon_{k,l}\omega)^p \wedge 
          \omega^{n-p} \]
and similarly for
$R_k^{(p)}$
and
$T_k^{(p)}$
one finally gets
\[ \int_{X-Sing(T_k^{(p)})} (T_k^{(p)} + \epsilon_{k,l}\omega)^p \wedge 
          \omega^{n-p} \leq
   \int_{X-Sing(T_{k,l})} (T_{k,l} + \epsilon_{k,l}\omega)^p \wedge 
          \omega^{n-p}. \]
Since
$\epsilon_{k,l} \rightarrow \epsilon_k$
the same line of arguments shows
\[ \int_{X-Sing(T_k^{(p)})} (T_k^{(p)} + \epsilon_{k,l}\omega)^p \wedge 
         \omega^{n-p}  \rightarrow 
   \int_{X-Sing(T_k^{(p)})} (T_k^{(p)} + \epsilon_k\omega)^p \wedge 
          \omega^{n-p}. \]
For 
$l$
big enough (depending on
$k$)
this gives
\[ s_k - \delta_k \leq 
   \int_{X-Sing(T_{k,l})} (T_{k,l} + \epsilon_{k,l}\omega)^p \wedge 
          \omega^{n-p} \leq s_{k+1}. \]

\noindent
Combining all these facts one gets a sequence of closed positive
$(1,1)-\!$
currents
$T_k = T_{k,l(k)}$
with analytic singularities in
$\alpha[-\epsilon_{k+1}\omega]$
such that the
$T_k$
converge weakly against
$T_{min}$,
and conditions (ii) and (iii) of the theorem are also satisfied.
\end{proof}

\noindent
\textbf{Remark.}
As long as
$T_{k,min} \rightarrow T_{min}$
weakly for
$k \rightarrow \infty$,
in the construction above it is not necessary that the
$T_{k,min}$
are computed w.r.t. the same smooth
$(1,1)-\!$
form on
$\alpha$.

\smallskip

\noindent
The approximation may be used e.g. to prove
\begin{lem} \label{limLelong-lem}
Let
$X$
be a compact K\"ahler manifold and
$\alpha \in H^{1,1}(X,\mathbb{R})$
a pseudo-effective class. Let
$\Delta^n \cong U \subset X$
be an open subset, and let
$p: \Delta^n \rightarrow \Delta^{n-1}$
be the projection onto the last
$n-1$
coordinates. Then there is a pluripolar set
$E \subset \Delta^{n-1}$
such that for all fibers
$\Delta$
over points in
$\Delta^{n-1} \setminus E$
\[ \lim_{\epsilon \downarrow 0} \inf_T \nu(T_{|\Delta},x) = \nu(\alpha,x)\ \ 
   \mathrm{for\ all\ } x \in \Delta, \]
where the 
$T$'s
run through all currents in
$\alpha[-\epsilon\omega]$
with analytic singularities, for which the restriction to
$\Delta$
is well-defined.
\end{lem} 
\begin{proof}
The proof is an application of the theory of
$(L,h)$-
resp.
$T$-
general curves. If
$T$
is an almost psoitive
$(1,1)$-
current on
$X$,
a smooth curve
$C$
(compact or not) will be called 
$T$-
general iff the restriction of
$T$
on
$C$
is well-defined and
\begin{itemize}
\item[(i)]
$C$
intersects no codim-2-component in any of the Lelong number level sets
$E_c(T)$,
\item[(ii)]
$C$
intersects every prime divisor 
$D \subset E_c(T)$
in the regular locus
$D_{reg}$
of this divisor,
$C$
does not intersect the intersection of two such prime divisors, and every 
intersection point
$x$
has the minimal Lelong number 
$\nu(T,x) = \nu(T,D) := \min_{z \in D} \nu(T,z)$,
\item[(iii)]
for all
$x \in C$,
the Lelong numbers
\[ \nu(T_{|\Delta},x) = \nu(T,x). \]
\end{itemize}

\noindent
Then theorem 2.1. in~\cite{E02} states that in a family of curves over
a smooth base there is a pluripolar subset in the base such that every curve
over points outside this pluripolar set is 
$T$-
general. In particular, this is true for currents
$T_k$
approximating
$T_{min}$
as in the theorem above. Since the union of countably many pluripolar sets is 
again pluripolar, this proves the lemma. 
\end{proof}

\subsection{Upper bound for the numerical dimension}

\noindent
The numerically trivial foliation w.r.t. a pseudo-effective class may be also
used to bound its numerical dimension, provided that the singularities of the
foliation are nice enough:
\begin{thm} \label{AppBd-thm}
Let
$X$
be a compact K\"ahler manifold with K\"ahler form 
$\omega$
and
$\alpha \in H^{1,1}(X,\mathbb{R})$
a pseudo-effective class. Let
$\mathcal{F}$
be the numerically trivial foliation w.r.t. 
$\alpha$
and suppose that the singularities of
$\mathcal{F}$
are isolated points. Then the numerical dimension
$\nu(\alpha)$ 
is less or equal to the codimension of the leaves of 
$\mathcal{F}$.
\end{thm}
\begin{proof}
Applying theorem~\ref{NumdimApp-thm} ,
one gets a sequence of closed
$(1,1)$-
currents with analytic singularities in
$\alpha[-\epsilon_k\omega]$
such that
\[ \lim_{k \rightarrow \infty} \int_{X - \mathrm{Sing\ }T_k} 
   (T_k + \epsilon_k\omega)^p \wedge \omega^{n-p} = 
   (\alpha^p.\omega^{n-p})_{\geq 0} \]
for all
$p = 1, \ldots, n$.
In these integrals,
the 
$T_k$'s
may be replaced by the residue currents 
\[ R_k = T_k - \sum \nu(T_k,D)[D] \]
of the Siu decomposition of the
$T_k$.

\noindent
Now the proof consists of two steps: first, let
$\Delta^n \cong U \subset X$
be an open set such that
$p: U \cong \Delta^n \rightarrow \Delta^l$
describes the numerical trivial foliation w.r.t. 
$\alpha$
locally in
$U$. 
Then use as in proposition~\ref{BdLKL-prop} that the 
$R_k$'s
get close to pulled back currents from the base
$\Delta^l$
to show

\smallskip

\noindent
\textit{Claim 1.}
For
$l < p \leq n$
and an open subset 
$U^\prime \subset\subset U$,
\[ \int_{U^\prime} (R_k + \epsilon_k\omega)^p \wedge \omega^{n-p} \rightarrow
   0. \]
\begin{proof}
Every
$R_k + \epsilon_k\omega$
may be written as a sum
$\sum_{i,j} \theta_{ij}^k dz_i \wedge d\overline{z}_j$.
Then every coefficient of
$(R_k + \epsilon_k\omega)^p$
w.r.t. the base
$dz_I \wedge d\overline{z}_J$
(with multi-index notation) is a product of
$p$
of these
$\theta_{ij}^k$.
If
$p > l$,
then one of these
$\theta_{ij}^k$
has index
$i \leq n-l$
or
$j \leq n-l$.

\noindent
As in proposition~\ref{BdLKL-prop} one can argue with the Schwarz inequality 
that
\[ |\theta_{ij}^k| \leq |\theta_{ii}^k|^{\frac{1}{2}} \cdot 
   |\theta_{jj}^k|^{\frac{1}{2}}. \]
Furthermore, let
$F_i$
be a sufficiently general fiber of the projection
$\Delta^n \rightarrow \Delta^{n-1}$
onto all but the ith coordinate,
$i = 1, \ldots , n-l$.
Since 
$R_k$
is a current with analytic singularities only in codimension
$2$,
a sufficiently general 
$F_i$
does not hit the singularities of
$R_k$.
Then
$\theta_{ii|F_i \cap U^\prime}^k$
is smooth and positive, and the numerical triviality implies that
\[ \int_{F_i \cap U^\prime} |\theta_{ii}^k| dz_i \wedge d\overline{z}_i =
   \int_{F_i \cap U^\prime} \theta_{ii}^k dz_i \wedge d\overline{z}_i
   \stackrel{k \rightarrow \infty}{\longrightarrow} 0. \]
Hence using Fubini for the ith or jth coordinate and applying the Schwarz 
inequality again, one gets that all the integrals of the terms of
$(R_k + \epsilon_k\omega)^p \wedge \omega^{n-p}$
tend to 
$0$
for
$k \rightarrow \infty$.
\end{proof}

\noindent
The second step is to give an estimate of the considered integrals around the
isolated singularities 
of the foliation by using the uniform boundedness of the Lelong numbers of 
(almost) positive currents in the same cohomology class.

\smallskip

\noindent
\textit{Claim 2.}
There is a sequence of compact sets
$K_i \subset X$
exhausting
$X - \mathrm{Sing\ }\mathcal{F}$
and a constant 
$C > 0$
such that for all
$1 \leq p \leq n$
\[ \int_{X - K_i} (R_k + \epsilon_k\omega)^p \wedge \omega^{n-p} \leq
   \delta_i, \]
and 
$\lim_{i \rightarrow \infty} \delta_i = 0$.
\begin{proof}
This is just an expanded version of Boucksom's argument 
in~\cite[Lem 3.1.11]{Bou02}. Choose a finite covering of
$X$
by open charts
$U_i$
isomorphic to the unity ball
$B \subset \mathbb{C}^n$,
such that the balls with half of the diameter still cover
$X$.
If 
$z^{(i)}$
denote coordinates on
$U_i$
one may find two constants
$C_1, C_2 > 0$
such that
\[ C_1\omega \leq \frac{i}{2}\partial\overline{\partial} |z^{(i)}|^2 \leq 
   C_2\omega \]
in
$V_i$,
for all
$i$.

\noindent
If
$x \in X$
lies in
$V_i$,
the Lelong number 
$\nu((R_k + \epsilon\omega)^p,x)$
is by definition the decreasing limit for
$r \rightarrow 0$
of
\[ \nu((R_k + \epsilon\omega)^p,x,r) := \frac{1}{(\pi r^2)^{n-p}} 
   \int_{|z^{(i)}-x|<r} (R_k + \epsilon\omega)^p \wedge 
   (\frac{i}{2}\partial\overline{\partial} |z^{(i)}|^2)^{p}. \]
On the one hand, for
$r \leq r_0$
one has
\[ \nu((R_k + \epsilon\omega)^p,x,r) \leq \nu((R_k + \epsilon\omega)^p,x,r_0) 
   \leq \frac{C_2}{(\pi r_0^2)^{n-p}} \int_X 
   (R_k + \epsilon\omega)^p \wedge \omega^{n-p}. \]
But
$\int_X (R_k + \epsilon\omega)^p \wedge \omega^{n-p} \leq 
 \int_X (T_k + \epsilon\omega)^p \wedge \omega^{n-p}$,
and the last integral depends only on the cohomology class of
$T_k$,
since
$\omega$
is closed.

\noindent
On the other hand,
\[ (\pi r^2)^{n-p} \nu((T + \epsilon\omega)^p,x,r) \geq 
   C_1 \int_{|z^{(i)}-x|<r} (T + \epsilon\omega)^p \wedge \omega^{n-p}. \]
For
$p < n$
the claim follows since 
$\mathrm{Sing\ }\mathcal{F}$
is compact, hence consists of only finitely many points. For
$p < n$
there is nothing to argue, since
$\nu(\alpha) = n$
implies that
$\alpha$
is big (\cite[Thm. 3.1.31]{Bou02}). Hence the numerically trivial foliation
coincides with the Iitaka fibration w.r.t.
$\alpha$,
because it is the identity map.
\end{proof}

\noindent
Both claims together show the theorem.
\end{proof}

\section{Surface Examples} \label{Ex-CEx-sec}

\noindent
If one constructs the numerical trivial foliation w.r.t. an incomplete system 
of currents with analytic singularities in
$\alpha[-\epsilon_k\omega]$,
$\epsilon_k \rightarrow 0$, 
then the leaf dimension is greater or equal than that of the numerical trivial
foliation w.r.t.
$\alpha$.
Unfortunately, the author could not prove any criterion when the leaf dimension
remains the same (hence the two foliations are equal). In general,
it seems quite difficult to decide whether a given foliation is numerically
trivial w.r.t. some pseudo-effective class
$\alpha$.
In the first two surface examples which follow, some ad-hoc arguments are used 
to show the identity of the constructed foliations and the numerically trivial
foliation w.r.t. the given pseudo-effective classes.

\subsection{A nef line bundle without smooth positive curvature form}
\label{firstex-ssec}

\noindent
This example was already discussed in~\cite{DPS94}: Let
$\Gamma = \mathbb{C}/(\mathbb{Z}+\mathbb{Z}\tau)$,
$\mathrm{Im}\tau > 0$, 
be an elliptic curve and let
$E$
be the rank 2 vector bundle over
$\Gamma$
defined by
\[ E = \mathbb{C} \times \mathbb{C}^2/(\mathbb{Z}+\mathbb{Z}\tau) \]
where the action is given by the two automorphisms
\[ \begin{array}{l}
   g_1(x,z_1,z_2) = (x+1,z_1,z_2)\\
   g_\tau(x,z_1,z_2) = (x+\tau,z_1+z_2,z_2),
   \end{array}\] 
and where the projection
$E \rightarrow \Gamma$
is induced by the first projection
$(x,z_1,z_2) \mapsto x$.
Then 
$\mathbb{C} \times \mathbb{C} \times \{0\}/(\mathbb{Z}+\mathbb{Z}\tau)$ 
is a trivial line subbundle
$\mathcal{O} \hookrightarrow E$,
and the quotient
$E/\mathcal{O} \cong \Gamma \times \{0\} \times \mathbb{C}$
is also trivial. Let
$L$
be the line bundle
$L = \mathcal{O}_E(1)$
over the ruled surface
$X = \mathbb{P}(E)$. 
The exact sequence
\[ 0 \rightarrow \mathcal{O} \rightarrow E \rightarrow \mathcal{O} 
     \rightarrow 0 \]
shows that
$L$
is nef over
$X$.

\noindent
Now, in \cite{DPS94} all hermitian metrics 
$h$
(including singular metrics) are determined such that the curvature current
$\Theta_h(L)$
is semipositive (in the sense of currents): These metrics 
have all the same curvature current
\[ \Theta_h(L) = [C], \]
where 
$C$
is the curve on
$X$
induced by 
$\{z_2 = 0\}$.
(This implies in particular that there exists no \textbf{smooth} positive 
hermitian metric on
$L$.)
To exclude the possibility that there exist positive currents in
$c_1(L)$
which are not the curvature current of a metric on
$L$
one proves the following
\begin{lem} \label{metricnumtriv-lem}
Let
$X$
be a projective complex manifold and 
$L$
a holomorphic line bundle on
$X$. 
Then for every closed positive current in
$c_1(L)$
there is a possibly singular hermitian metric 
$h$
on
$L$
such that the curvature current
\[ \Theta_h(L) = T. \]
\end{lem} 
\begin{proof}
Let
$T$
be any positive current in
$c_1(L)$.
By \cite{Bon95} there exists a line bundle 
$L^\prime$
on
$X$
with a possibly singular hermitian metric 
$h^\prime$
such that
$\Theta_{h^\prime}(L^\prime) = T$.
(This is just the usual construction of a cycle in 
$H^1(X, \mathcal{O}^\ast)$). 
The line bundle
$N = (L^\prime)^{-1} \otimes L$
is numerically trivial, hence nef. Consequently there exists a positive 
singular hemitian metric 
$h_N$
on
$N$
such that the class of the curvature current
\[ \{\Theta_{h_N}(N)\} = 0 \in H^{1,1}(X, \mathbb{R}). \]

\noindent
Now, all closed positive currents in 
$0 \in H^{1,1}(X, \mathbb{R})$
have the form
$dd^c \phi$
for some plurisubharmonic function on
$X$.
Since
$\phi$
is upper semi-continuous it attains its supremum. But then the maximum 
principle implies that
$\phi$
is a constant function. Therefore the only closed positive current in 
$0 \in H^{1,1}(X, \mathbb{R})$
is the zero form. This implies
$\Theta_{h_N}(N) = 0$
(as a current).

\noindent
Furthermore this gives the hermitian metric
$h = h_N \otimes h^\prime$
on
$L = N \otimes L^\prime$
with
$\Theta_h(L) = T$.
\end{proof}

\noindent
So 
$[C]$
really is a positive current with minimal singularities in
$c_1(L)$. 
But then
$X$
is numerically trivial w.r.t. 
$[C]$,
and the associated numerical trivial foliation has only one leaf
$X$
with codimension
$0$.

\noindent
On the other hand, 
$L$
is certainly not numerically trivial since it intersects a fiber of
$X = \mathbb{P}(E)$
with intersection number
$1$.
Consequently, the moving intersection number
$(c_1(L) \cdot c_1(L))_{\geq 0} = c_1(L) \cdot c_1(L)$
is strictly positive,
and 
$(X,c_1(L))$
is a counter example to equality of the numerically trivial foliation w.r.t.
the positive closed 
$(1,1)$-
current with minimal singularities and that w.r.t. the associated 
pseudo-effective cohomology class.

\noindent
Now there is an obvious candidate for a numerically trivial foliation w.r.t.
$c_1(L)$: 
its leaves are the projection of the curves
$\mathbb{C} \times \{p\}$
in
$\mathbb{P}_C(E)$. 
The strategy to show this has two parts: first, one constructs a sequence of 
currents 
$T_k \in c_1(L)[-\epsilon_k\omega]$
for some K\"ahler form
$\omega$
on
$X$
and a sequence
$\epsilon_k$
of positive real numbers tending to
$0$
such that the foliation mentioned above is the numerically trivial foliation 
w.r.t. this sequence of
$T_k$'s. 
Second, one uses that the restriction of the 
$T_k$'s
to any
$\mathbb{P}^1$-\!
fiber of
$\mathbb{P}_C(E)$
is 
$\geq c \cdot \omega$,
for some fixed number
$c > 0$.

\noindent
The construction of the 
$T_k$
requires a careful study of almost positive (singular) hermitian metrics 
$h$
on
$L$:
As the total space of
$L^{-1}$
is equal to
$E^\ast$
blown up along the zero section, the function
\[ \phi(\zeta) = \log \parallel \zeta \parallel^2_{h^{-1}},\ \zeta \in L^{-1} 
   \]
associated to any hermitian metric
$h$
on
$L$
can also be seen as a function on
$E^\ast$
satisfying the log-homogeneity condition
\[ \phi(\lambda\zeta) = \log |\lambda| + \phi(\zeta)\ \mathrm{for\ every\ }
   \lambda \in \mathbb{C}. \]
One has
\[ \frac{i}{2\pi} \partial\overline{\partial} \phi(\zeta) = 
   \pi^\ast_{L^{-1}}\Theta_h(L),\ \pi_{L^{-1}}: L^{-1} \rightarrow X. \]
Thus
$\Theta_h(L)$
is almost positive iff 
$\phi$
is almost plurisubharmonic on
$E^\ast$.

\noindent
The total space of
$E^\ast$
is the quotient
$E^\ast = \mathbb{C} \times \mathbb{C}^2/(\mathbb{Z} + \mathbb{Z}\tau)$
by the dual action
\[ \begin{array}{rcl}
   g_1^\ast(x,w_1,w_2) & = & (x+1, w_1,w_2) \\
   g_\tau^\ast(x,w_1,w_2) & = & (x+\tau, w_1,w_1+w_2).
   \end{array} \]
The function
$\phi$
gives rise to a function
$\widetilde{\phi}$
on 
$\mathbb{C} \times \mathbb{C}^2$
which is invariant by
$g_1^\ast, g_\tau^\ast$
and log-homogeneous w.r.t. 
$(w_1,w_2)$,
and
$\widetilde{\phi}$
is almost plurisubharmonic iff
$\phi$
is almost psh. Even more is true: Interpret
$X$
as the zero section of the total space of
$L^{-1}$
and let
$\omega_X, \omega_{L^{-1}}$
be positive
$(1,1)$-\!
forms on
$X, L^{-1}$.
Then there are constants
$C_1, C_2 > 0$
such that
\[ 0 \leq \pi^\ast_{L^{-1}} \omega_X \leq C_1 \omega_{L^{-1}},\ 
   0 \leq \omega_{L^{-1}|X} \leq C_2 \omega_X. \]
Hence 
$\pi^\ast_{L^{-1}} \Theta_h \geq -\epsilon \omega_{L^{-1}}$
implies
$\Theta_h \geq -\epsilon C_2 \omega_X$,
and
$\Theta_h \geq -\epsilon\omega_X$
implies
$\pi^\ast_{L^{-1}} \Theta_h \geq -\epsilon C_1 \omega_{L^{-1}}$.
Consequently, instead of constructing currents
$T_k \geq -\epsilon_k\omega_X$,
$\epsilon_k \rightarrow 0$
on
$X$,
it suffices to construct currents
$\Theta_k \geq -\epsilon_k^\prime \omega_{L^{-1}}$,
$\epsilon_k^\prime \rightarrow 0$,
and functions
$\widetilde{\phi}_k$
on
$\mathbb{C} \times \mathbb{C}^2$
such that
$i\partial\overline{\partial}\widetilde{\phi}_k = \Theta_k$
and the 
$\widetilde{\phi}_k$ 
are invariant by
$g_1^\ast, g_\tau^\ast$
and log-homogeneous w.r.t. 
$(w_1,w_2)$.

\noindent
This is done by using a
gluing procedure developed in \cite{Dem92}:
\begin{lem} \label{glue-lem}
Let 
$U_j^\prime \subset\subset U_j^{\prime\prime} \subset\subset  U_j$
be locally finite open coverings of a complex manifold 
$X$
by relatively compact open sets, and let
$\theta_j$
be smooth nonnegative functions with support in
$U_j^{\prime\prime}$
such that
$\theta_j \leq 1$
on
$U_j^{\prime\prime}$
and
$\theta_j = 1$
on
$U_j^\prime$.
Let
$A_j \geq 0$
be such that
\[ i(\theta_j\partial\overline{\partial}\theta_j - 
     \partial\theta_j \wedge \overline{\partial}\theta_j) \geq - A_j\omega\ 
   \ \ \ \ \mathrm{on\ } U_j^{\prime\prime}\setminus U_j^\prime. \] 
Finally, let
$w_j$
be almost psh functions on
$U_j$
with the property that
$i\partial\overline{\partial}w_j \geq \gamma$
for some real 
$(1,1)-\!$
form
$\gamma$
on
$X$,
and let
$C_j$
be constants such that
\[ w_j(x) \leq C_j + \sup_{k \neq j, U_k^\prime \ni x} w_k(x)\ \ \ \ \ 
  \mathrm{on\ } U_j^{\prime\prime}\setminus U_j^\prime. \]
Then the function
$w = \log(\sum \theta_j^2e^{w_j})$
is almost psh and satisfies
\[ i\partial\overline{\partial}w \geq \gamma - 
   2\left( \sum_j \mathbf{1}_{U_j^{\prime\prime}\setminus U_j^\prime}
                  A_j e^{C_j}\right) \omega. \]
\end{lem}
\begin{proof}
See \cite[Lemma 3.5]{Dem92}. The proof is reproduced here, because the results
of the following computation are needed later on: Setting
$\alpha_j = \theta_j\partial w_j + 2\partial \theta_j$,
one gets
\[ \begin{array}{rcl}
   \partial w & = & 
   \frac{\sum \theta_j e^{w_j} \alpha_j}{\sum \theta_j^2 e^{w_j}} \\
   \partial\overline{\partial}w & = & 
   \frac{\sum (\alpha_j \wedge \overline{\alpha}_j +
               \theta_j^2\partial\overline{\partial}w_j + 
               2\theta_j\partial\overline{\partial}\theta_j -
               2\partial\theta_j \wedge \overline{\partial}\theta_j)e^{w_j}}
        {\sum \theta_j^2 e^{w_j}} - 
   \frac{\sum_{j,k} \theta_j e^{w_j} \theta_k e^{w_k} 
         \alpha_j \wedge \overline{\alpha}_k}{(\sum \theta_j^2 e^{w_j})^2}  \\
    & = & \frac{\sum_{j < k} (\theta_j\alpha_k-\theta_k\alpha_j) \wedge 
     (\theta_j\overline{\alpha}_k-\theta_k\overline{\alpha}_j)e^{w_j}e^{w_k}}
               {(\sum \theta_j^2 e^{w_j})^2} + 
   \frac{\sum \theta_j^2 e^{w_j}\partial\overline{\partial}w_j}
        {\sum \theta_j^2 e^{w_j}} + \\
    &  & + \frac{\sum (2\theta_j\partial\overline{\partial}\theta_j -
                 2\partial\theta_j \wedge \overline{\partial}\theta_j)e^{w_j}}
                {\sum \theta_j^2 e^{w_j}}. \\
   \end{array} \]
The first term is
$\geq 0$,
the second term is
$\geq \gamma$.
If a point 
$x$
lies in the support of the third term, it is contained in some 
$U_j^{\prime\prime}\setminus U_j^\prime$. 
By assumption there is a 
$k$
such that
$w_j(x) \leq C_j + w_k(x)$
and
$x \in U_k^\prime$,
hence
$\theta_k(x) = 1$.
This implies 
\[ \frac{\sum (2\theta_j\partial\overline{\partial}\theta_j -
                 2\partial\theta_j \wedge \overline{\partial}\theta_j)e^{w_j}}
                {\sum \theta_j^2 e^{w_j}} \geq -2 \sum 
   \mathbf{1}_{U_j^{\prime\prime}\setminus U_j^\prime} e^{C_j} A_j \omega. \] 
\end{proof}

\noindent
On any open subset of
$\Gamma$
on which the
$\mathbb{P}^1$-\!
bundle is trivial, one may define a metric on
$L$
as the pullback w.r.t. a trivialization of the Fubini-Study metric on
$\mathbb{P}^1$.
These metrics induce the candidates for the 
$w_j$:
Set
\[ \begin{array}{rcl}
   U_j & = & \{(x,w_1,w_2): j\mathrm{Im}\tau < \mathrm{Im}x < 
                        (j+1)\mathrm{Im}\tau \}, \\
   U_j^\prime & = & \{(j+\frac{1}{6})\mathrm{Im}\tau < 
                      \mathrm{Im}x < (j+\frac{3}{4})\mathrm{Im}\tau,\ 
                      |jw_1+w_2| > |(j+1)w_1+w_2| \}, \\
   U_j^{\prime\prime} & = & \{(j+\frac{1}{12})\mathrm{Im}\tau < 
                 \mathrm{Im}x < (j+\frac{5}{6})\mathrm{Im}\tau,\ 
                 (1+\epsilon)|jw_1+w_2| > |(j+1)w_1+w_2| \} 
   \end{array} \]
and take
\[ \widetilde{u}_j = \frac{1}{2} \log (|w_1|^2 + |jw_1+w_2|^2), \]
set
\[ \begin{array}{rcl}
   V_j & = & \{(x,w_1,w_2): j\mathrm{Im}\tau < \mathrm{Im}x < 
                        (j+1)\mathrm{Im}\tau \}, \\
   V_j^\prime & = & \{(j+\frac{1}{4})\mathrm{Im}\tau < 
                      \mathrm{Im}x < (j+\frac{5}{6})\mathrm{Im}\tau,\ 
                      |jw_1+w_2| < |(j+1)w_1+w_2| \}, \\
   V_j^{\prime\prime} & = & \{(j+\frac{1}{6})\mathrm{Im}\tau < 
                 \mathrm{Im}x < (j+\frac{11}{12})\mathrm{Im}\tau,\ 
                 |jw_1+w_2| < (1+\epsilon)|(j+1)w_1+w_2| \} 
   \end{array} \]
and take
\[ \widetilde{v}_j = \frac{1}{2} \log (|w_1|^2 + |(j+1)w_1+w_2|^2), \]
set
\[ \begin{array}{rcl}
   W_j & = & \{(x,w_1,w_2): (j-\frac{1}{2})\mathrm{Im}\tau < \mathrm{Im}x < 
                        (j+\frac{1}{2})\mathrm{Im}\tau \}, \\
   W_j^\prime & = & \{(j-\frac{1}{4})\mathrm{Im}\tau < 
                      \mathrm{Im}x < (j+\frac{1}{4})\mathrm{Im}\tau\}, \\
   W_j^{\prime\prime} & = & \{(j-\frac{1}{3})\mathrm{Im}\tau < 
                 \mathrm{Im}x < (j+\frac{1}{3})\mathrm{Im}\tau\} 
   \end{array} \]
and take
\[ \widetilde{w}_j = \frac{1}{2} \log (|w_1|^2 + |jw_1+w_2|^2). \]
A straight forward computation shows that these sets and functions satisfy all
the conditions of the glueing lemma. Choose smooth nonnegative
functions 
$\theta_j, \theta_j^\prime, \theta_j^{\prime\prime}$ 
invariant under multiplication of
$(w_1, w_2)$
by
$\lambda \in \mathbb{C}$ 
such that
$g_1^\ast \theta_j = \theta_j$
and
$g_\tau^\ast \theta_j = \theta_{j+1}$
(similarly for
$\theta_j^\prime, \theta_j^{\prime\prime}$)
and constants
$A_j, A_j^\prime , A_j^{\prime\prime} > 0$
satisfying the conditions of the glueing lemma for some positive 
$(1,1)$-\!
form
$\widetilde{\omega}$
on
$\mathbb{C} \times \mathbb{C}^2$.
Applying the glueing lemma one gets an almost plurisubharmonic function
\[ \widetilde{w} = \log (\sum \theta_j^2 e^{\widetilde{u}_j} +
                         \sum (\theta_j^\prime)^2 e^{\widetilde{v}_j} +
                     \sum (\theta_j^{\prime\prime})^2 e^{\widetilde{w}_j}), \]
which is invariant by
$g_1^\ast, g_\tau^\ast$
and log-homogeneous.

\noindent
This construction also works for 
$k\widetilde{u}_j,k\widetilde{v}_j,k\widetilde{w}_j$
and then gives 
$\widetilde{w}^{(k)}$.
Setting
\[ \widetilde{\phi}_k := \frac{1}{k} \widetilde{w}^{(k)}, \]
one gets almost plurisubharmonic functions
$\widetilde{\phi}_k$
with negative part arbitrary small (compared to a positive 
$(1,1)$-\!
form
$\widetilde{\omega}$
as above).

\noindent
Since the
$\widetilde{u}_j,\widetilde{v}_j,\widetilde{w}_j$
do not depend on
$x$,
one has
\[ \partial_x \widetilde{u}_j = \partial_x \widetilde{v}_j =
   \partial_x \widetilde{w}_j = 0. \]
Furthermore there is a constant
$B > 0$
such that
\[ 2\theta_j\partial\overline{\partial}\theta_j -
                 2\partial\theta_j \wedge \overline{\partial}\theta_j  < 
   B \omega,
   \]
similarly for
$\theta_j^\prime, \theta_j^{\prime\prime}$,
and a constant
$C > 0$
such that
\[ 4|\theta_j \partial \theta_{j^\prime} - 
     \theta_{j^\prime} \partial \theta_j|^2  < C \omega \]
and all other possible pairs of
$\theta_j, \theta_j^\prime, \theta_j^{\prime\prime}$.
Consequently, in
$x$-\!
direction the first and the third term of
$i\partial\overline{\partial}\widetilde{\phi}_k$
in the glueing lemma are bounded from above by
$\frac{3C}{k}\omega, \frac{B}{k}\omega$,
by the same trick as in the proof of the glueing lemma (the 
$3$
comes from the fact that every
$x$
lies in the support of at most
$\theta_j, \theta_j^\prime, \theta_j^{\prime\prime}$
for exactly one
$j$),
and the second term vanishes. Hence the foliation under consideration is really
the numerically trivial foliation w.r.t. the currents induced on
$X$.  

\noindent
On the other hand, it follows from the construction that the restriction of 
the induced currents
$T_k$
to the
$\mathbb{P}^1$-\!
fibers of
$X = \mathbb{P}(E)$
remain
$> \epsilon\omega$
for some
$\epsilon > 0$.

\noindent
Let
$T_k^\prime \in \alpha[-\epsilon_k\omega]$
be another sequence of currents representing
$\alpha$.
If
$\Delta^2 \cong U \subset X$
is an open subset with coordinates 
$z_1, z_2$
such that the lines 
$\{z_1 = a\}$
belong to
$\mathbb{P}^1$-\!
fibers and
$\{z_2 = b\}$
are subsets of the leaves of the foliation one can write 
\[ T_k + \epsilon_k\omega = \sum_{i,j = 1}^2 \theta_{ij}^{(k)} 
                           idz_i \wedge d\overline{z}_j,\ 
   T_k^\prime  + \epsilon_k\omega = \sum_{i,j = 1}^2 \theta_{ij}^{\prime(k)} 
                                    idz_i \wedge d\overline{z}_j. \]
By the remark above,
\[ (\theta_{22}^{(k)})_{|\{z_1 = a\}} idz_2 \wedge d\overline{z}_2 > 
   \epsilon\omega \]
for all
$a$,
and
\[ \widetilde{\theta}^{(k)} := \theta_{11}^{(k)} idz_1 \wedge d\overline{z}_1 +
               \theta_{12}^{(k)} idz_1 \wedge d\overline{z}_2 +
               \theta_{21}^{(k)} idz_2 \wedge d\overline{z}_1  
   \stackrel{k \rightarrow \infty}{\longrightarrow} 0 \]
by the numerical triviality (use as before the Schwarz inequality for the 
terms 
with~$\theta_{12}^{(k)}, \theta_{21}^{(k)}$).

\noindent
Since the numerical dimension of 
$L$
is
$1$,
one knows furthermore that
\[ \lim_{k \rightarrow \infty} \int_{X-\mathrm{Sing\ }T_k^\prime} 
   (T_k + \epsilon_k\omega) \wedge (T_k^\prime + \epsilon_k\omega) = 0. \]
But 
\[ (T_k + \epsilon_k\omega) \wedge (T_k^\prime + \epsilon_k\omega) = 
   \widetilde{\theta}^{(k)} \wedge (T_k^\prime  + \epsilon_k\omega) + 
   \theta_{22}^{(k)}idz_2 \wedge d\overline{z}_2 \wedge 
   \theta_{11}^{\prime(k)} idz_1 \wedge d\overline{z}_1, \]
hence the vanishing of the limits above implies
\[ \int_{(\Delta^\prime)^2-\mathrm{Sing\ }T_k^\prime} \theta_{11}^{\prime(k)} 
   idz_1 \wedge d\overline{z}_1 \wedge idz_2 \wedge d\overline{z}_2 
   \stackrel{k \rightarrow \infty}{\longrightarrow} 0, \]
where
$\Delta^\prime \subset\subset \Delta$
is any open disk such that
$(\Delta^\prime)^2 \subset U \cong \Delta^2$.

\noindent
Consequently,
$\int_{\Delta_b^\prime-\mathrm{Sing\ } T_k^\prime} 
 (T_k^\prime + \epsilon_k\omega) 
 \stackrel{k \rightarrow \infty}{\longrightarrow} 0$
for almost all
$b \in \Delta^\prime$
(where 
$\Delta_b^\prime = \{b\} \times \Delta^\prime$).
The definition of the numerically trivial foliation requires that 
$\int_{\Delta_b^\prime-\mathrm{Sing\ }T_k^\prime}
 (T_k^\prime + \epsilon_k\omega) 
 \stackrel{k \rightarrow \infty}{\longrightarrow} 0$
for \textbf{all} 
$b \in \Delta^\prime$.
To prove this one can use the same line of arguments as in the proof of the 
Local Key Lemma for pseudo-effective classes: One tries to show that
\[ \lim_{k \rightarrow \infty}   
   | \int_{\Delta_b^\prime-\mathrm{Sing\ }T_k^\prime}
     (T_k^\prime + \epsilon_k\omega) -
     \int_{\Delta_0^\prime-\mathrm{Sing\ }T_k^\prime}
     (T_k^\prime + \epsilon_k\omega) | = 0. \]
Following the proof of proposition~\ref{BdLKL-prop} one sees that it is enough
to show that
\[ \lim_{k \rightarrow \infty} \int_{\partial\Delta^\prime \times\Delta_{0,b}} 
   |\theta_{11}^{\prime (k)}|dV \cdot 
   \int_{\partial\Delta^\prime \times \Delta_{0,b}} 
   |\theta_{22}^{\prime (k)}|dV = 0, \]
where
$\Delta_{0,b}$
is the disk with center in
$b/2$
and radius
$|b/2|$,
and
$dV$
is a volume element of
$\partial\Delta^\prime \times\Delta_{0,b}$.

\noindent
As in the proof of proposition~\ref{BdLKL-prop} there is a bound
$M > 0$
such that for all
$k$
there is a disk 
$\Delta_k^\prime \subset\subset \Delta$
containing
$\Delta^\prime$
with
\[ \int_{\partial\Delta_k^\prime \times \Delta^\prime} 
   |\theta_{22}^{\prime (k)}|dV < M. \]
For the first term, look at the 
$(1,1)$-\!
form
$\eta = idz_2 \wedge d\overline{z}_2$
and take a disk
$\Delta^\prime \subset\subset \Delta^{\prime\prime} \subset\subset \Delta$.
Then by the arguments above,
\[ \int_{(\Delta^{\prime\prime}-\Delta^\prime) \times \Delta^\prime} 
            |\theta_{11}^k| dV =
   \int_{(\Delta^{\prime\prime}-\Delta^\prime) \times \Delta^\prime} 
   (T_k^\prime + \epsilon_k\omega) \wedge \eta 
   \stackrel{k \rightarrow \infty}{\longrightarrow} 0. \]

\noindent
By Fubini, one gets a disk
$\Delta_k^\prime$
such that
\[ \int_{\partial\Delta_k^\prime \times \Delta^\prime} 
            |\theta_{11}^k| dV 
   \stackrel{k \rightarrow \infty}{\longrightarrow} 0, \]
and one concludes that the limit above is indeed
$0$.

\subsection{Mumford's example} \label{MumEx-ssec}

\noindent
Back to our counter example at the beginning: it is easy to construct a 
closed positive
$(1,1)-\!$
current on
$L = \mathcal{O}_{\mathbb{P}^1}(1)$
such that the leaves of the associated numerically trivial foliation are 
1-dimensional. Take a measure
$\omega$
invariant w.r.t. the representation of
$\pi(C)$
in
$\mathrm{PGL}(2)$.
This gives a measure on
$(\Delta \times \mathbb{P}^1)/\pi(C)$
transversal to the foliation induced by the images of
$\Delta \times \left\{ p \right\}$. 
Averaging out the integration currents of the leaves with this transverse 
measure gives an (even smooth) closed positive
$(1,1)-\!$
current in the first Chern class of
$L = \mathcal{O}_{\mathbb{P}(E)}(1)$
which vanishes on the leaves but not in any transverse direction.

\noindent
Most of this example is explained in the introduction; the only assertions 
not already discussed are the existence of a measure
$\omega$
in 
$c_1(\mathcal{O}_{\mathbb{P}^1}(1))$
invariant w.r.t. the unitary representation of
$\pi(C)$
in
$\mathrm{GL}(2)$
and the smoothness of the metric which results from averaging out the 
integration currents of the leaves. But this is easy, too: Take the Haar 
measure 
$\omega$
on the Lie group 
$U(2)$
which is absolutely continuous (\cite[Ch.14]{Dieu70}). Since
$U(2)$
operates transitively on
$\mathbb{P}^1$
this measure induces a
$U(2)-\!$ 
invariant measure on the homogeneous quotient space
$\mathbb{P}^1$. 
Since
$U(2)$
is compact it is possible to normalize
$\omega$
such that 
$\mathbb{P}^1$
has measure
$1$.
Hence averaging over the integration currents of the leaves w.r.t. 
$\omega$
gives a smooth positive
$(1,1)-\!$
form which is still in the first Chern class of
$L = \mathcal{O}_{\mathbb{P}(E)}(1)$.
Since it is smooth it is a current with minimal singularities on
$L$,
and obviously, this current is numerically trivial on the leaves.

\noindent
On the other hand it is strictly positive on the
$\mathbb{P}^1$-\!
fibers, hence the foliation is numerically trivial w.r.t. the cohomology class
by the same argument as in the first example. 

\noindent
\textbf{Remark.}
The difference to the previous example is that the unitary group is compact and
consequently its Haar measure is finite. This is not the case for the group of
linear automorphisms generated by 
$(z_1,z_2) \mapsto (z_1 + z_2,z_2)$.

\subsection{$\mathbb{P}^2$ blown up in 9 points}

\noindent
Consider the following situation: Let
$C \subset \mathbb{P}^2$
be a smooth elliptic curve and let
$p_1, \ldots, p_8 \in C$
be sufficiently general points. The aim is to study the numerically trivial 
foliation w.r.t. the anticanonical bundle
$-K_X$
on varieties
$X_p = \mathbb{P}^2(p_1, \ldots, p_8,p)$
blown up in points
$p \in C$.

\noindent
Let
$E_i = \pi^{-1}(p_i)$
be the exceptional divisor on 
$X$
over
$p_i$.
First of all,
$-K_X = \mathcal{O}_{\mathbb{P}^2}(3) + \sum E_i$
is nef and
$-K_X^2 = 0$.
Next, the pencil of elliptic curves on
$\mathbb{P}^2$
through
$p_1, \ldots, p_8$
has a base point
$q$.
So
$X_q =  \mathbb{P}^2(p_1, \ldots, p_8,q)$
is an elliptic fibration
$\pi_q: X_q \rightarrow \mathbb{P}^1$.
The pull back of a smooth positive metric on
$\mathcal{O}_{\mathbb{P}^1}(1)$
gives a smooth semipositive hermitian metric on
$-K_{X_q}$
which is strictly positive in directions transverse to the fibers. Hence by the
same arguments as in the two examples above, the fibration is the
numerically trivial foliation w.r.t.
$-K_{X_q}$.  

\noindent
For points 
$p \neq q$
in
$C$
there is only one section
in
$-K_{X_p}$,
the strict transform
$C^\prime$
of
$C$. 
But if one considers torsion points (w.r.t. to
$q$)
of order 
$m$
on
$C$
then a calculation in \cite{DPS96} shows that
$-mK_{X_p}$
defines again an elliptic fibration over
$\mathbb{P}^1$.
This fibration yields a smooth semipositive hermitian metric on
$-mK_{X_p}$,
hence on
$-K_{X_p}$,
and again the fibration is the numerically trivial foliation w.r.t.
$-mK_{X_p}$.

\noindent
The question is: What happens if non-torsion points
$p \in C$
are blown up ? In particular: Is there always a smooth semipositive hermitian
metric on
$-K_{X_p}$
inducing a holomorphic foliation on
$X_p$,
which may be seen as the limit of the fibrations of
$X_{p_k}$ 
where the
$p_k$
are torsion points ? (The last question was asked in \cite{DPS96}.) 
A strategy to answer it is to 
use the theory of holomorphic foliations on surfaces, as developed e.g. in 
\cite{Bru00}. 
\begin{Def}
A (holomorphic) \textit{foliation}
$\mathcal{F}$
on a compact complex surface
$X$
is a coherent analytic rank 1 subsheaf
$T_{\mathcal{F}}$
of the tangent bundle
$T_X$
(the tangent bundle of the foliation)
fitting into an exact sequence
\[ 0 \rightarrow T_{\mathcal{F}} \rightarrow T_X \rightarrow 
     \mathcal{J}_Z \otimes N_{\mathcal{F}} \rightarrow 0 \]
for a suitable invertible sheaf
$N_{\mathcal{F}}$
(the normal bundle of the foliation)
and an ideal sheaf
$\mathcal{J}_Z$
whose zero locus consists of isolated points called the \textit{singularities}
$\mathrm{Sing}(\mathcal{F})$
of
$\mathcal{F}$.
\end{Def}

\noindent
Furthermore, one can easily show that
$T_{\mathcal{F}}^\ast \otimes N_{\mathcal{F}}^\ast = K_X$.

\noindent
Numerically trivial foliations 
$\{\mathcal{F},(U_i,p_i)\}$
on surfaces
$X$
with
$\mathcal{F}$
of rank 1 are such foliations: If
$\mathcal{F}$
is not a line bundle then replace it by
$\mathcal{F}^{\ast\ast}$. 
As a reflexive sheaf on a surface this is a line bundle \cite[1.1.10]{OSS}, 
and dualizing the inclusion
$\mathcal{F} \subset T_X$
twice shows that it is still a subsheaf of
$T_X$.
Furthermore, 
$\mathcal{F}$
is locally integrable because it has rank 1, hence the maps 
$p_i$
exist trivially. 

\noindent
Let
$\mathcal{X}$
be
$\mathbb{P}^2(p_1, \ldots, p_8) \times C$
blown up in the diagonal
\[ \Delta_{C \times C} \subset C \times C \subset 
   \mathbb{P}^2(p_1, \ldots, p_8) \times C. \]
The fibers of
$\mathcal{X}$
over 
$p \in C$
are just the
$X_p$
for all
$p$.
If there is an algebraic family of foliations on the
$X_p$
such that over torsion points, the foliation coincides with the fibration 
described above, then (at least generically) the conormal line bundles 
$N_{F_p}^\ast$
should also fit into a family. But this is impossible, as the following 
computation shows:
\begin{lem}
Let
$C, q, X_p$
be as above, and let 
$p$
be a torsion point w.r.t. 
$q$
of order
$m$.
Let
$N_{F_p}$
be the normal bundle of the foliation induced by the fibration
$\pi_p: X_p \rightarrow \mathbb{P}^1$.
Then
\[ N_{F_p}^\ast \cong (m+1)K_{X_p}. \] 
\end{lem}
\begin{proof}
Let
$D$
be an irreducible component of a fiber of
$\pi = \pi_p$
with multiplicity 
$l_D$.
If
$\eta$
is a local non-vanishing 
$1$-
form on
$\mathbb{P}^1$
then
$\pi^\ast(\eta)$
is a local section of
$\pi^\ast(K_{\mathbb{P}^1})$
vanishing of order
$l_D-1$
on
$D$.
Hence,
\[ N_{F_p}^\ast = \pi^\ast(K_{\mathbb{P}^1}) \otimes 
                  \mathcal{O}_{X_p}(\sum (l_D-1)D). \]
The relative canonical bundle formula (for elliptic fibrations, 
see~\cite{Fri98}) tells that
\[ K_{X_p} = \pi^\ast(K_{\mathbb{P}^1} \otimes 
             (R^1_\ast\pi\mathcal{O}_{X_p})^\ast) \times 
             \mathcal{O}_{X_p}(\sum (l_F-1)F), \]
where the sum is taken over all fibers 
$F$
occuring with multiplicity
$l_F$
in the fibration. 

\noindent
There are two differences between the two formulas: First, in the relative 
canonical bundle formula occurs the term
\[ L := (R^1_\ast\pi\mathcal{O}_{X_p})^\ast. \]
Now, 
$\deg L \geq 0$,
and
$\deg L = 0$
would imply that 
$L$
is a torsion bundle on 
$\mathbb{P}^1$,
hence it is trivial, and 
$X_p = C \times \mathbb{P}^1$ --
a contradiction.
If
$L$
is nontrivial, a short calculation with spectral sequences shows that
\[ 0 = p_g = \deg L - g(\mathbb{P}^1) + 1, \]
hence 
$\deg L = 1$,
and
$L = \mathcal{O}_{\mathbb{P}^1}(1)$
(see~again~\cite[Ch.VII]{Fri98}). This shows
\[ \pi^\ast(K_{\mathbb{P}^1} \otimes L) = 
   \pi^\ast \mathcal{O}_{\mathbb{P}^1}(-1) = mK_{X_p}, \]
and together with the relative canonical bundle formula this shows that
$mC$
is the only multiple fiber.

\noindent
The second difference is that some fibers may contain multiple components, but
are not multiple themselves. By the classification of singular fibers of 
elliptic fibrations this is only possible if there are
$-2$-\!
curves (\cite{Fri98}). But on
$\mathbb{P}^2$
blown up in 
$9$
points in general position, there are no
$-2$-\!
curves. Hence
\[ \mathcal{O}_{X_p}(\sum (l_D-1)D) = \mathcal{O}_{X_p}(\sum (l_F-1)F),\]
and the claim of the lemma follows.
\end{proof}

\noindent
The threefold
$\mathcal{X}$
is also a counter example to equality of numerical dimension and codimension 
of the leaves of the numerically trivial foliation w.r.t. some 
pseudo-effective class: Set
\[ L := -\pi^\ast(p_1^\ast K_{\mathbb{P}^2(p_1, \ldots, p_8)}) - E_\Delta +
        p_2^\ast \mathcal{O}(nr), \]
where
$p_1$
is the projection of
$\mathbb{P}^2(p_1, \ldots, p_8) \times C$
onto
$\mathbb{P}^2(p_1, \ldots, p_8)$,
$p_2$
is the projection of
$\mathcal{X}$
onto
$C$,
$r$
is any point on
$C$
and
$n > 0$
an integer. The restriction of 
$L$
to any fiber over 
$p \in C$
is the anticanonical bundle
$K_{X_p}^\ast$.

\noindent
For 
$n$
sufficiently big,
$L$
is nef:
$L$
is effective, since 
$D = C \times C + nX_r$
is contained in
$|L|$.
Consequently, to prove the nefness of 
$L$
it suffices to show that all curves
$E \subset C \times C$
have non-negative intersection number with
$L$.
To this purpose first get an overview over all curves on
$C \times C$:
According to the general theory of abelian surfaces the Picard number of
$C \times C$
is
$4$
or
$3$
depending on whether 
$C$
has complex multiplication or not (\cite[2.7]{CT}. Hence it suffices to look
at the fibers of the two projections of
$C \times C$
onto
$C$,
the diagonal, and if necessary, on some other curve constructed as the graph of
complex multiplication in
$C \times C$. 
Since it is a graph of an isomorphism, such a curve maps isomorphically to
$C$
under both projections.

\noindent
Now, one has to compute the degree of the restriction of
$L$
to
$E$.
This restriction may also be seen as the restriction of the divisor
$D_{|D}$
to such an
$E$.
Let
$C^\prime$
be a sufficiently general curve in the pencil
$|-K_{\mathbb{P}^2(p_1, \ldots, p_8)}|$.
Then the strict transform of
$C^\prime \times C$
is an element of
$-\pi^\ast(p_1^\ast K_{\mathbb{P}^2(p_1, \ldots, p_8)})$
and intersects
$C \times C$
in
$\{q\} \times C$. 
Furthermore,
$E_\Delta$
intersects
$C \times C$
in the diagonal
$\Delta_{C \times C}$. 
Therefore,
\[ D_{|D} \sim \{q\} \times C + n(C^\prime \times \{r\}) - 
                     \Delta_{C \times C} + n(C \times \{r\}), \]
where
$E_r$
is the exceptional divisor over
$r$
in
$X_r$.
And
$L$
is nef if 
$n$ 
is 
$\geq$
the maximum of
$1$
(this is the intersection number of fibers
$C \times \{p\}$
with the diagonal) and the intersection number of the curve coming from complex
multiplication (if existing) with the diagonal. (The self intersection number 
of the diagonal is
$0$
since the tangent bundles on
$C \cong \Delta_{C \times C}$
and
$C \times C$
are trivial.)  
\begin{prop}
Let
$\mathcal{X}, L$
be as above. Then the numerical dimension 
$\nu(L)$
of
$L$
is 
$2$,
but the numerically trivial foliation w.r.t. 
$c_1(L)$
is the identy map.
\end{prop}
\begin{proof}
To prove
$L^2 \neq 0$,
observe that
$L^2$
is represented by the cycles in the expression above for
$D_{|D}$.  
This is not
$\equiv 0$,
since the intersection number with
$\{q\} \times C$
is positive for
$n \geq 1$.

\noindent
The numerically trivial
foliation w.r.t. 
$c_1(L)$
cannot be the trivial map onto a point, because in fibers 
$X_p$
over torsion points
$p$
there are curves which are not numerically trivial. Since immersed disks which 
do not lie in a fiber of the projection onto 
$C$ 
are not numerically trivial, the only possible numerically trivial
foliation w.r.t. 
$c_1(L)$
with
$2$-\!
dimensional leaves is the fibration onto
$C$.
But this is impossible by the same reason as above. 
To exclude the possibility that the numerically trivial foliation has
$1$-\!
dimensional leaves, one notes first that over torsion points
$p$,
the fibers of
$\pi_p: X_p \rightarrow \mathbb{P}^1$
are numerically trivial: This is clear since these fibers 
$F$
are projective, hence
$\int_F T_k$
only depends on the cohomology class of the 
$T_k$,
and
$\int_F c_1(L)$
is certainly 
$0$.

\noindent
This can be used to show that the
$1$-\!
dimensional leaves of a numerically trivial foliation must lie in the fibers
$X_p$
of
$\mathcal{X}$:
Otherwise, let 
$\Delta^3 \cong U \subset \mathcal{X}$
be any open subset with coordinates 
$x,z_1,z_2$
such that the projection onto
$C$
is given by the projection onto the first coordinate, and the foliation is 
described by the projection onto the two last coordinates. Choose
$x$
such that 
$x = 0$
corresponds to a torsion point 
$p_0$. 
Shrinking
$U$
if necessary, one can suppose that the fibers of 
$\pi_{p_0}$
are smooth in
$U$.
But then the Local Key Lemma for pseudo-effective classes implies that there
are
$2$-\!
dimensional numerically trivial leaves, contradiction.

\noindent
Next one shows that the 
$1$-\!
dimensional leaves in fibers
$X_p$,
where
$p$
is a torsion point, must be the fibers of
$\pi_p: X_p \rightarrow \mathbb{P}^1$:
Take an ample line bundle 
$A$
on
$\mathcal{X}$. 
Since
$L$
is nef,
$L^k \otimes A$
is also ample, and some multiple is very ample. The global sections of this 
very ample line bundle generate a smooth metric on
$L^k \otimes A$
whose strictly positive curvature form may be written as
$k(T_k + \frac{1}{k}\omega_A)$,
for some form 
$T_k \in c_1(L)[-\frac{1}{k}\omega_A]$.

\noindent
Let
$p \in C$
be any torsion point of order 
$m$
and
$\pi_p: X_p \rightarrow \mathbb{P}^1$
the induced fibration. Let
$T = i\partial\overline{\partial}\log(|z_1|^2+|z_2|^2)$
be a strictly positive curvature form in
$c_1(\mathcal{O}_\mathbb{P}^1(1))$.
Then
\[ (T_k + \frac{1}{k}\omega_A)_{|X_p} \geq \frac{1}{m} \pi_p^\ast T. \]
But this means in particular that for any disk
$\Delta \subset X_p$
not immersed into a fiber of
$\pi_p$,
\[ \int_\Delta T_k + \frac{1}{k}\omega_A \geq 
   \frac{1}{m} \int_\Delta \pi_p^\ast T > 0. \]
Hence the leaves of the numerically trivial foliation w.r.t.
$c_1(L)$
coincide with the fibers of
$\pi_p$
in
$X_p$.

\noindent
But this is impossible, as shown above.
\end{proof}

\noindent
\textbf{Remark.}
This proposition does not exclude the possibility that (some of)
the
$X_p$
over non-torsion points
$p$
have a numerically trivial foliation with 1-dimensional leaves. 

\noindent
Another result dealing with this type of foliations is
\begin{prop}[Brunella]
Let
$\mathcal{F}$
be a foliation on a compact algebraic surface
$X$
and suppose that
$\mathcal{F}$
is tangent to a smooth elliptic curve
$E$,
free of singularities of
$\mathcal{F}$.
Then either
$E$
is a (multiple) fiber of an elliptic fibration or, up to ramified coverings 
and birational maps, 
$\mathcal{F}$
is the suspension of a representation
$\rho: \pi_1(\widehat{E}) \rightarrow \mathrm{Aut}(\mathbb{CP}^1)$,
$\widehat{E}$
an elliptic curve.
\end{prop}

\appendix

\section{Singular foliations}

\noindent
One can define foliations on complex manifolds as 
involutive subbundles of the tangent bundle. Then the classical theorem of
Frobenius asserts that through any point there is a unique integral 
complex submanifold~\cite{Miy86}. Singular foliations may be defined as 
involutive coherent subsheaves of the tangent bundle, which are furthermore 
saturated, that is, their quotient with the tangent bundle is torsion free.
In the points where the rank is maximal, one may use again the Frobenius 
theorem to get leaves. Since in this paper the reasoning is always explicitely
using the leaves their existence is directly incorporated in the definition of
a singular foliation:
\begin{Def}\label{singfol-def}
Let
$X$
be an
$n-\!$
dimensional compact complex manifold.
Let
$\mathcal{F} \subset T_X$
be a saturated subsheaf of the tangent bundle with maximal rank
$k$
and
$Z \subset X$
be the analytic subset where 
$\mathcal{F}/m_{X,x}\mathcal{F} \rightarrow T_{X,x}$
is not injective.

\noindent
$\mathcal{F}$
induces a singular foliation described by the following data:
$X-Z$
is covered by open sets
$U_i \cong \Delta^n$
such that for the smooth holomorphic map
$p_i: U_i \rightarrow \Delta^{n-k}$
coming from the projection
$\Delta^n \rightarrow \Delta^{n-k}$,
\[ \mathcal{F}_{|U_i} = T_{U_i/\Delta^{n-k}}. \]
Such a foliation will be denoted by
$\{\mathcal{F},(U_i,p_i)\}$.
\end{Def}

\noindent
Next, one defines the inclusion relation for numerically trivial foliations as
above
\begin{Def}
A numerically trivial foliation is contained in another one,
\[ \left\{ \mathcal{F}, (U_i, p_i) \right\} \sqsubset
   \left\{ \mathcal{G}, (V_j, q_j) \right\}, \]
iff there is a Zariski open set
$U \in X$
such that
$\mathcal{F}_{|U} \subset \mathcal{G}_{|U}$.
\end{Def}

\noindent
In particular this means that the leaves of
$\left\{ \mathcal{F}, (U_i, p_i) \right\}$
are contained in those of
$\left\{ \mathcal{G}, (V_j, q_j) \right\}$.

\noindent
The next aim is to construct a common refinement 
$\left\{ \mathcal{H}, (W_k, r_k: W_k \rightarrow \Delta^{n-m}) \right\}$
of two singular foliations
$\left\{ \mathcal{F}, (U_i, p_i: U_i \rightarrow \Delta^{n-k}) \right\}$,
$\left\{ \mathcal{G}, (V_j, q_j: V_j \rightarrow \Delta^{n-l}) \right\}$,
that is
\[ \left\{ \mathcal{F}, (U_i, p_i) \right\} \sqsubset 
   \left\{ \mathcal{H}, (W_k, r_k) \right\},\ 
   \left\{ \mathcal{G}, (V_j, q_j) \right\} \sqsubset
   \left\{ \mathcal{H}, (W_k, r_k) \right\}. \]
To this purpose one has first to analyze the local 
picture when two foliations meet transversally everywhere: Let
$W$
be a complex manifold with two isomorphisms
$q_1: W \rightarrow \Delta^n$,
$q_2: W \rightarrow \Delta^n$.
Let
$p_1: W \rightarrow \Delta^{n-k}$,
$p_2: W \rightarrow \Delta^{n-l}$
be the composition of 
$q_1$,
$q_2$
with the projections of
$\Delta^n$
onto the last 
$n-k$
resp.
$n-l$
factors. 

\noindent
\hspace{4cm}
\xymatrix{
   W \ar[r]^{p_1} \ar[d]_{p_2} & \Delta^{n-k} \\
   \Delta^{n-l}
         }

\vspace{0.5cm}

\noindent
Suppose that the 
$p_1-\!$
and
$p_2-\!$
fibers intersect transversally everywhere. 

\noindent
If
$n \leq k + l$,
by choosing appropriate coordinates 
$p_1$
will be the projection on the first
$n-k$
coordinates while
$p_2$
is the projection on the last
$n-l$
coordinates. In particular, the smallest fibration
$p$
whose fibers contain all fibers of
$p_1$
and
$p_2$
is the trivial fibration onto a point.

\noindent
So suppose from now on that
$n > k + l$.
Again by choosing appropriate coordinates via the implicit function theorem 
and possibly further restricting 
$W$
one can describe the configuration of the two sets of fibers in the following 
way (look at the next figure): The horizontal sections of 
$p_1$
consist of 
$p_2-\!$
fibers which are parallel hyperplanes, and over each point
$y \in \Delta^{n-k}$
the
$p_2-\!$
fibers through the points in
$q^{-1}(y)$
project into a pencil of hyperplanes through a common central hyperplane
$\Delta^{k^\prime} \subset \Delta^{n-k}$
containing
$y$.

\noindent
This central hyperplane is isomorphic to
$\Delta^{k^\prime}$
for all
$y \in \Delta^{n-k}$,
and different central hyperplanes are parallel in
$\Delta^k$.
Let
$r: \Delta^{n-k} \rightarrow \Delta^{n - k - k^\prime}$
be the projection with the central hyperplanes as fibers.
Consequently the new projection
\[ p: W \stackrel{p_1}{\rightarrow} \Delta^{n-k} 
        \stackrel{r}{\rightarrow} \Delta^{n - k - k^\prime} \] 
is the smallest fibration whose fibers contain both the fibers of
$p_1$
and
$p_2$.

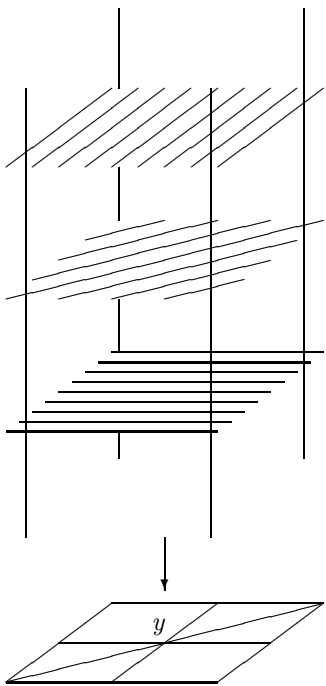
\begin{center}
\begin{picture}(130,260)(0,0)
\put(7.5,55){\line(0,1){170}}
\put(77.5,55){\line(0,1){170}}
\put(42.5,85){\line(0,1){10}}
\put(42.5,125){\line(0,1){20}}
\put(42.5,175){\line(0,1){20}}
\put(42.5,225){\line(0,1){30}}
\put(112.5,85){\line(0,1){170}}

\put(0,95){\line(1,0){80}}
\put(5,98.7){\line(1,0){80}}
\put(10,102.5){\line(1,0){80}}
\put(15,106.2){\line(1,0){80}}
\put(20,110){\line(1,0){80}}
\put(25,113.7){\line(1,0){80}}
\put(30,117.5){\line(1,0){80}}
\put(35,121.2){\line(1,0){80}}
\put(40,125){\line(1,0){80}}

\put(0,195){\line(4,3){40}}
\put(10,195){\line(4,3){40}}
\put(20,195){\line(4,3){40}}
\put(30,195){\line(4,3){40}}
\put(40,195){\line(4,3){40}}
\put(50,195){\line(4,3){40}}
\put(60,195){\line(4,3){40}}
\put(70,195){\line(4,3){40}}
\put(80,195){\line(4,3){40}}

\put(0,145){\line(4,1){120}}
\put(20,145){\line(4,1){90}}
\put(40,145){\line(4,1){60}}
\put(60,145){\line(4,1){30}}
\put(100,175){\line(-4,-1){90}}
\put(80,175){\line(-4,-1){60}}
\put(60,175){\line(-4,-1){30}}

\put(60,55){\vector(0,-1){20}}

\put(0,0){\line(1,0){80}}
\put(40,30){\line(1,0){80}}
\put(0,0){\line(4,3){40}}
\put(80,0){\line(4,3){40}}

\put(55,20){$\mathit{y}$}

\put(20,15){\line(1,0){80}}
\put(0,0){\line(4,1){120}}
\put(40,0){\line(4,3){40}}
 
\end{picture}
\end{center}

\vspace{0.5cm}

\noindent
So outside the singularities of
$\mathcal{F}$
and
$\mathcal{G}$
and the locus where the leaves of the foliation do not intersect 
transversally or even coincide, it is clear how to define the common 
refinement. 

\noindent
To get a better feeling  for the locus of the other points, look at the 
following two-dimensional toy example, where the two foliations are marked with
dotted and dashed lines.

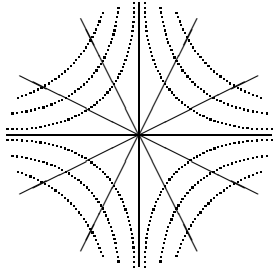
\begin{center}
\begin{picture}(100,100)(-50,-50)

\put(0,50){\line(0,-1){100}}

\put(-50,0){\line(1,0){100}}

\put(0,0){\line(1,2){22}}
\put(0,0){\line(2,1){45}}

\put(0,0){\line(1,-2){22}}
\put(0,0){\line(2,-1){45}}

\put(0,0){\line(-1,2){22}}
\put(0,0){\line(-2,1){45}}

\put(0,0){\line(-1,-2){22}}
\put(0,0){\line(-2,-1){45}}

\qbezier[50](2,50)(2,2)(50,2)
\qbezier[40](8,48)(12,12)(48,8)
\qbezier[30](14,46)(22,22)(46,14)

\qbezier[50](2,-50)(2,-2)(50,-2)
\qbezier[40](8,-48)(12,-12)(48,-8)
\qbezier[30](14,-46)(22,-22)(46,-14)

\qbezier[50](-2,50)(-2,2)(-50,2)
\qbezier[40](-8,48)(-12,12)(-48,8)
\qbezier[30](-14,46)(-22,22)(-46,14)

\qbezier[50](-2,-50)(-2,-2)(-50,-2)
\qbezier[40](-8,-48)(-12,-12)(-48,-8)
\qbezier[30](-14,-46)(-22,-22)(-46,-14)

\end{picture}
\end{center}

\noindent
At least, the exceptional points form an analytic subset of
$X$: 
A point
$x \in U_i \cap V_j$
is contained in this set iff the differential of
$p_i \times q_j: U_i \cap V_j \rightarrow \Delta^{n-k} \times \Delta^{n-l}$
in
$x$ 
has not full rank, that is iff all maximal minors of this differential vanish
in
$x$. 
But it still remains the task to define a saturated subsheaf 
$\mathcal{H} \subset T_X$
which locally coincides with the relative tangential sheaf of the projections
$r_k: W_k \rightarrow \Delta^{n-m}$. 

\noindent
To do this one goes back to the purely algebraic definition of (singular)
foliations: The subsheaf
$\mathcal{F} \subset T_X$
induces such a foliation iff it is involutive, that is, closed under the Lie
bracket, which means
$[\mathcal{F},\mathcal{F}] \subset \mathcal{F}$.
Then there is a natural candidate for a subsheaf defining the union of the 
foliations given by
$\mathcal{F}$
and
$\mathcal{G}$:
the smallest involutive subsheaf
$\mathcal{H} \subset T_X$
containing both
$\mathcal{F}$
and
$\mathcal{G}$. 
It exists because it may be constructed as the subsheaf generated by
$\mathcal{F}$,
$\mathcal{G}$
$[\mathcal{F},\mathcal{G}]$,
$[[\mathcal{F},\mathcal{G}],\mathcal{F}]$,
$[[\mathcal{F},\mathcal{G}],\mathcal{G}]$,
$[[\mathcal{F},\mathcal{G}],[\mathcal{F},\mathcal{G}]]$
and so on.
\begin{lem}
Let
$\left\{ \mathcal{F}, (U_i, p_i: U_i \rightarrow \Delta^{n-k}) \right\}$,
$\left\{ \mathcal{G}, (V_j, q_j: V_j \rightarrow \Delta^{n-l}) \right\}$
be two singular foliations on an
$n$-\!
dimensional complex manifold
$X$.
Let
$x$
be a point not in the analytic subset
$Z \subset X$
consisting of the singular locus of
$\mathcal{F}$,
$\mathcal{G}$
and the points where the leaves of
$\mathcal{F}$
and
$\mathcal{G}$
do not intersect transversally. Then on the common refinement
$r_k: W_k \rightarrow \Delta^{n-m}$
around
$x$
constructed as above, the smallest involutive subsheaf
$\mathcal{H}$
containing both
$\mathcal{F}$
and
$\mathcal{G}$
coincides with the relative tangential sheaf of
$r_k$.
\end{lem}
\begin{proof}
Since the two foliations intersect transversally around
$x$,
it is obvious that the smallest saturated involutive subsheaf in
$T_{X|W_k}$
containing both
$\mathcal{F}_{|W_k}$
and
$\mathcal{G}_{|W_k}$
is the relative tangent sheaf of the projection 
$r_k$.
Glueing together one gets a saturated involutive subsheaf
$\mathcal{H}_{X-Z} \subset T_{X-Z}$
on the open set
$X-Z$.

\noindent
Let
$U \cong \Delta^n$
be a neighborhood of some point
$z \in Z$,
and let
$H = \{f = 0\}$,
$f \in \mathcal{O}(U)$,
be an analytic hyperplane in
$U$
containing the analytic subset 
$Z \cap U$.
Now,
$\mathcal{O}(U-H) = \mathcal{O}(U)_f$,
and one can define the sections of
$\mathcal{H}$
on
$U$
as the intersection
\[ \mathcal{H}(U) = T_X(U) \cap \mathcal{H}_{X-Z}(U-H) \]
in
$T_X(U-H) \cong \mathcal{O}(U)^n$.
Since
$T_X(U)$
and
$\mathcal{H}_{X-Z}(U-H)$
are involutive,
$\mathcal{H}(U) \subset T_X(U)$
is closed under the Lie bracket, too. Furthermore,
$\mathcal{H}(U)$
is the smallest \textbf{saturated} submodule of
$T_X(U)$
such that
$\mathcal{H}(U)_f = \mathcal{H}_{X-Z}(U-H)$,
as the following algebraic lemma shows. 

\noindent
Finally, since the same is true for
$\mathcal{F}(U)$
and
$\mathcal{G}(U)$,
they are both contained in
$\mathcal{H}(U)$.
\end{proof}

\begin{lem}
Let
$R$
be a commutaive integral ring,
$f \in R$,
and
$M^f \subset R_f^k$
a submodule such that
$R_f^k/M^f$
is torsion free. Then
$M = M^f \cap R^k$
is the smallest submodule of
$R^k$
such that
$M_f = M^f$
and
$R^k/M$
is torsion free.
\end{lem}
\begin{proof}
If
$N \subset M^f \subset R$
such that there exists
$m \in M - N$,
but still
$N_f = M^f$,
then
$m \in N_f$.
Hence there is an
$n \in N$
and
$l \in \mathbb{N}$
such that
$m = \frac{n}{f^l}$
or
$m \cdot f^l = n$. 
But then
$f$
is a torsion element of
$R^k/N$.
\end{proof}

\noindent
This shows that
$\left\{ \mathcal{H}, (W_k, r_k: W_k \rightarrow \Delta^{n-m}) \right\}$
is really a singular foliation.

\newcommand{\etalchar}[1]{$^{#1}$}

\end{document}